\documentclass[letterpaper,11pt]{article}
\usepackage{verbatim}
\usepackage{times}
\usepackage{empheq}
\usepackage{amsmath, amssymb}
\usepackage{amsthm}
\usepackage[T1]{fontenc}
\usepackage[english]{babel}
\usepackage[utf8]{inputenc}
\usepackage{graphicx}
\usepackage{hyperref}
\usepackage{mathrsfs}
\usepackage{mathabx}
\usepackage{color}
\usepackage{ulem}
\usepackage{anysize}
\usepackage{subfigure}
\usepackage{natbib}

\marginsize{2.5cm}{2.5cm}{2.5cm}{2.5cm}
\newtheorem{Theo}{Theorem}

\newtheorem{Prop}{Proposition}
\newtheorem{Lemma}{Lemma}[section]
\newtheorem{Corollary}{Corollary}
\theoremstyle{definition}
\newtheorem{Example}{Example}[section]
\newtheorem{Rem}{Remark}
\usepackage{upgreek}

\usepackage{setspace}
\usepackage{physics}
\AtBeginDocument{}
\addto\captionsfrench{}
\DeclareMathOperator\sgn{sgn}
\DeclareMathOperator\supp{supp}
\usepackage[toc,page]{appendix}

\title{A general framework for SPDE-based stationary random fields}

\author
{\\ \ \\ 
	 Ricardo Carrizo Vergara$^1$, Denis Allard$^2$, Nicolas Desassis$^1$\\
     \normalsize{$^{1}$Geostatistics, MINES ParisTech, 77350, Fontainebleau, France.}\\
	\normalsize{$^{2}$Biostatistics and Spatial Processes, BioSP, INRA, 84914, Avignon, France.}\\
}
\date{\today}
\begin{document}
	\maketitle

\noindent {\bf Abstract} \quad This paper presents theoretical advances in the application of the Stochastic Partial Differential Equation (SPDE) approach in geostatistics. We show a general approach to construct stationary models related to a wide class of linear SPDEs, with applications to spatio-temporal models having non-trivial properties. Within the framework of Generalized Random Fields, a criterion for existence and uniqueness of stationary solutions for this class of SPDEs is proposed and proven. Their covariance are then obtained through their spectral measure. We present a result relating the covariance in the case of a White Noise source term with that of a generic case through convolution. Then, we obtain a variety of SPDE-based stationary random fields. In particular, well-known results regarding the Matérn Model and Markovian models are recovered. A new relationship between the Stein model and a particular
SPDE is obtained. New spatio-temporal models obtained from evolution SPDEs of arbitrary temporal derivative order are then obtained, for which properties of separability and symmetry can be controlled. We also obtain results concerning stationary solutions for physically inspired models, such as solutions to the heat equation, the advection-diffusion equation, some Langevin's equations and the wave equation. 

\bigskip


\newpage

\section{Introduction}

Finding new statistical models for analyzing spatio-temporal data that appropriately capture the complex interactions between space and time observed in natural phenomenon while allowing efficient computations able to handle very large datasets is a very active field of research. The typical approach in modeling a variable that varies spatio-temporally is to consider it as a realization  of a random field, i.e. a stochastic process indexed in space-time. The common practice is to describe its statistical properties by its covariance function which must be positive-definite, thereby limiting the choice of available models and making the construction of models with realistic features intricate.

Most of the commonly used space-time covariance models are built by modifying or combining generic covariance models defined for $\mathbb{R}^d$, $d=1,2,\ldots$.  These basic models are usually  stationary and isotropic. Commonly known generic models are covariance functions of exponential, powered exponential, Mat\'ern or Cauchy type, amongst many others \citep{ChilesDelfiner}. Space-time separable covariance models are constructed taking a tensor product between  a spatial and a temporal covariance. Separability is often an overly simplistic assumption, since it cannot capture sophisticated interaction between space and time. One of the first attempts to build non-separable covariance functions lead to the so-called product-sum class of models which simply adds and multiplies valid covariance models in the space and time domains  \citep{deIaco2001,deIaco2002,PorcuQuasi}. Even though this approach is perfectly valid from a mathematical standpoint, it is not grounded on physical considerations. The Gneiting class \citep{Gneiting2002} provides flexible non-separable  space-time models. Its construction is based on mixtures arguments, with no reference to physical considerations. Contrarily to the product-sum model, the Gneiting class implies higher space-time correlation than separability, in accordance to most observed phenomenon. 

Another important notion characterizing some spatio-temporal covariances is that of full symmetry \citep{Gneiting2006}. In a fully symmetric model, the direction of the time evolution is ignored, obtaining equal covariance values if we look either forward or backward in time. Separable covariance functions are necessarily fully symmetric, but not vice-versa. Product sum models and the Gneiting class are fully symmetric, non separable covariance functions. 
Atmospheric or environmental processes are often under the influence of prevailing air or water flows which are incompatible with full symmetry. Transport effects of this type can easily be modeled with the help of a purely spatial covariance function and a possibly random velocity vector. See \cite{Benoit2018} for an example of application to precipitation fields and \cite{Ailliot_etal2011} for an application on significant height wave fields with varying velocities. We refer to \cite{Gneiting2006} for a more detailed review on usual spatio-temporal covariance models.

The recent Stochastic Partial Differential Equation (SPDE) approach advocated in \cite{Lindgren2011} has open a new paradigm for handling large to very large ($> 10^6$) spatial datasets.  This approach consists in modeling the spatial variable as arising from the solution of a particular class of SPDEs for which the representation on a finite grid presents interesting Markov properties making it computable even for very large grids.  Specifically, following \cite{Whittle1963},  \cite{Lindgren2011} consider the following SPDE
\begin{equation}
	( \kappa^2 - \Delta)^{\frac{\alpha}{2}}U = W, 
	\label{eq:Lindgren2011}
\end{equation}
where $\alpha > \frac{d}{2}$, $\Delta$ is the Laplacian, $W$ is a White Noise process and $U$ is the unknown random field whose covariance function is a Matérn covariance. The SPDE approach is an important paradigm shift from both theoretical and practical perspective.   From a theoretical viewpoint, in contrast to the previous statistically oriented constructions reviewed above,  it proposes a physically grounded construction, for which the  parameters carry traditional physical interpretation such as diffusivity, reaction and transport \citep[see also][]{Whittle1963,Dong1990,Kelbert2005}. It allows the construction of models with interesting non-separability and non-symmetry properties  \citep{JonesAndZhang1997, Brown2000}. Non-stationarity can also easily be accounted for \citep{Fuglstad2013}. Other works based on this paradigm include \cite{Vecchia1985}, \cite{GayAndHeyde1990} and \cite{RuizMedina2016}. 
From a practical viewpoint, spatial prediction and simulation can be computed with methods brought from numerical analysis and PDE-solving methods such as the Finite Elements Method. The sparsity of the  matrices involved in the computations allows extremely efficient treatment of very large datasets for which classical geostatistical methods fail due to their high computational cost. We refer to \citet{Lindgren2011} for a detailed presentation of this framework and to \citet{Simpson2012} for a discussion on the advantages of this approach compared to classical techniques. Thanks to this unique combination of theoretical and practical properties, the SPDE approach has been widely used for analyzing large data sets, in particular in environment or climate science \citep{Bolin2011, Camaletti2013,Huang2017,Mena2017}.
This approach has also inspired the development of other PDE-solver based methods with efficient performances for a wider class of models \citep{Sigrist2015,Liu2016,Bolin2017}. \citet{Lang2011} propose an efficient method of simulation for solutions of some SPDEs based on Fourier Analysis, taking advantage of the low computational cost of the Fast Fourier Transform (FFT).  

Despite the huge potential of the SPDE approach, theoretical advances have been scarce in a spatio-temporal context, in particular because the requirement of a Markov structure for fast matrix calculations imposes some constraints. In R-INLA, which is the commonly used R package using the SPDE representation (\ref{eq:Lindgren2011}) of Gaussian fields, a temporal effect can be modeled as an autoregressive process or as a random walk; see \cite{Camaletti2013} for an application to particulate matters and \cite{Opitz2017} for a recent review on R-INLA with a focus on spatio-temporal applications.  \cite{Sigrist2015} show that the solution of a stochastic advection-diffusion partial differential equation provides a flexible class of models for spatio-temporal phenomena.  Inference and simulations of these models are computationally feasible for large data sets thanks to spatial FFT and filtering techniques in time. In a spatial context,  actual application of the SPDE approach relies mainly on Eq. (\ref{eq:Lindgren2011}) with its associated Matérn covariance. It is therefore of theoretical and practical interest to build, in a quite general setting, spatial and spatio-temporal covariance models that can be related to specific SPDEs, thereby offering a wide class of new theoretical models, which in some cases can provide  usual physical interpretability. In addition, the computational efficiency brought by numerical analysis solvers can be used to study these models. This setting should offer the possibility to handle a very general class of random fields and yet should be easy to use in order to simplify the conception, characterization and exploitation of new models. Clearly, it is expected that known results and known models will appear as special cases in this setting. 

For this purpose, this work proposes a general framework based on the theory of Generalized Random Fields \citep{Ito1954,Rozanov1982,Matheron1965} which is the stochastic analogue of Schwartz's theory of Distributions \citep{Schwartz1966}. This ancient theory, not so popular among the statistical community allows to rigorously study stationary models related to some SPDEs.
We express under which conditions a stationary covariance model corresponds to the unique stationary solution of some SPDE. When such a model exists, the associated spectral measure is exhibited. The well-known result by Whittle \citep{Whittle1963} is easily recovered in this framework. This framework also includes many other known models such as 
stationary Markov models \citep{Rozanov1977}, long-range dependent random fields \citep{GayAndHeyde1990,Anh1998} appearing as Matérn models without range parameter, the Stein spectral measure \citep{Stein2005} and solutions of some linear evolution equations \citep{Kelbert2005,Sigrist2015}. It also allows us to obtain new spatio-temporal models that can present non-trivial properties such as non-separability and non-symmetry. In particular, non-symmetric evolution models with fractional temporal regularity are obtained.

The rest of the paper is organized as follows. In section \ref{Section:MainResult} we present our main result which provides a criteria for the existence and uniqueness of stationary solutions for a wide class of linear SPDEs. In this Section, we use the minimum necessary concepts to rigorously present the result. Technical details are left to Section \ref{Section:TheoreticalFramework} where the complete theoretical framework is presented. It is based on the concept of Generalized Random Fields in a tempered framework. Here, the random field is no longer a function but a distribution \citep{Schwartz1966} for which differential operators and the Fourier transform are well-defined operations. The interest of this section is, in addition to the proof of our main result, the construction of a rigorous framework where operations on random fields are well-defined and relatively easy to use.  White Noise is then introduced in Section \ref{Section:WhiteNoiseSourceTerm}. We define it as a particular generalized random field playing a central role. We present an important result that relates the covariance of the solutions of the SPDEs with any source term to the covariance of the solution of the same SPDE with a White Noise source term. In sections \ref{Section:IllustrationsSpatial} and \ref{Section:SpatioTemporalModels} we show examples of models that can be conceived within our framework, involving both known and new models. In section \ref{Section:IllustrationsSpatial} we review some known cases in a spatial context, which involve the Matérn covariance and Markov models. In section \ref{Section:SpatioTemporalModels} we work in a spatio-temporal context. In section \ref{Section:SteinModel} we relate the Stein model \citep{Stein2005} to a particular SPDE
 and in section \ref{Section:EvolutionModels} we present a wide-class of new spatio-temporal stationary models arising as solutions to evolution equations and which present non-trivial properties. We show examples  having special interest both in physics and statistics. A general analysis of the properties based on the spectral measures of the models is made, in particular with regards to their spatial structure. We specify the corresponding spatial SPDE when possible.  We conclude in \ref{Section:Conclusion} with some final words.

\section{Presentation of the main theoretical result}
\label{Section:MainResult}

In this section we introduce our main theoretical result regarding the existence and the uniqueness of stationary solutions for a rich class of linear SPDEs. All models presented later in this work, either in a spatial context in Section \ref{Section:IllustrationsSpatial}, or in a spatio-temporal one in Section \ref{Section:SpatioTemporalModels}, derive from this construction. This offers a unified framework to a variety of spatial and spatio-temporal models that have been presented or revisited recently. This result is presented here in general, not completely rigorous terms. A more detailed presentation is voluntarily deferred to Section \ref{Section:TheoreticalFramework},  where all proofs and formal definitions are given.

\subsection{Introduction}
\label{Section:IntroductionMainResult}

A second order stationary real random function over $\mathbb{R}^{d}$ is a family of squared-integrable real random variables indexed over the euclidean space, $Z = (Z(x))_{x \in \mathbb{R}^{d}}$, such that its mean function $m_{Z}(x) = \mathbb{E}( Z(x) )$ is constant and its covariance function $C_{Z}(x,y) = \mathbb{C}ov( Z(x) , Z(y) ) $ depends only on the gap $x-y$.  Without loss of generality, we will consider  that $m_{Z}=0$. The function $\rho_{Z} : \mathbb{R}^{d} \to \mathbb{R}$ such that $\rho_{Z}(x-y) = \mathbb{C}ov( Z(x),Z(y) )$ is called the \textit{stationary covariance function} and it is positive-definite. By Bochner's Theorem \citep[see for example ][chapter 37]{Donoghue1969}, $\rho_{Z}$ is the Fourier transform of a positive, finite and even measure over $\mathbb{R}^{d}$, $\mu_{Z}$, referred to as the \textit{spectral measure} of $Z$: $\rho_{Z} = \mathscr{F}(\mu_{Z})$, where $\mathscr{F}$ denotes the Fourier transform on $\mathbb{R}^{d}$ (see Appendix \ref{Section:TemperedDistribution} for the convention of the Fourier transform used in this work). The covariance function $\rho_{Z}$ and the spectral measure $\mu_{Z}$ can equivalently be used to fully characterize the covariance structure of $Z$. 

In this work, it will be necessary to consider more general mathematical objects that allow us to deal rigourously with linear differential operators and Fourier transforms on random fields. We will use \textit{Generalized Random Fields} (GeRF), which are an analogous to the generalization of functions presented in Schwartz's theory of Distributions; see for example \citet{Ito1954} for a theory of stationary GeRFs. In this framework, the random fields have only meaning when applied to test functions in some particular functional space, and not necessarily when evaluated in points of the space. We present all the technical details in section \ref{Section:TheoreticalFramework}. For now, we mention that the covariance structure of a stationary GeRF can be described by not necessarily finite spectral measures. To characterize those, we consider the class  ${\cal M}_{SG}^+(\mathbb{R}^d)$ of slow-growing positive (Borel) measures over $\mathbb{R}^d$. Members of  ${\cal M}_{SG}^+(\mathbb{R}^d)$ can have infinite total mass, but they grow at most at a polynomial rate. Specifically, 
\begin{equation}
\label{eq:SlowGMeasure}
{\cal M}_{SG}^+(\mathbb{R}^d):= \left\{ \mu \hbox{ positive measure over } \mathbb{R}^{d} \ \big| \quad \int_{\mathbb{R}^{d}}(1+|x|^{2})^{-N}d\mu(x) < \infty \quad \hbox{for some }  N \in \mathbb{N} \right\}.
\end{equation}
If $\mu_{Z} \in \mathcal{M}_{SG}^{+}(\mathbb{R}^{d})$ is even, it can be used as a spectral measure of a real stationary GeRF $Z$, and its Fourier transform (in distributional sense) $\rho_{Z}=\mathscr{F}(\mu_{Z})$ is called the \textit{stationary covariance distribution} of $Z$ (we will omit the adjective \textit{stationary} when it is clear from context). This distribution is not necessarily a continuous function and thus $Z$ is not necessarily a random function with pointwise meaning. However, when the condition in Eq. (\ref{eq:SlowGMeasure}) is satisfied for $N = 0$, we are back to the usual framework of second order stationary random functions and to Bochner's characterization of continuous stationary covariance functions. From now on, every even measure in $\mathcal{M}^{+}_{SG}(\mathbb{R}^{d})$ will be said to be a spectral measure.  

The focus in this work is on a quite general class of linear stochastic equations which encompasses those considered in \cite{Whittle1963} \cite{Lindgren2011}, \cite{Sigrist2015}, \cite{Bolin2011}, \cite{Anh1998}, \cite{Kelbert2005} and \cite{GayAndHeyde1990}.  Linear operators involved here are not strictly speaking  \textit{differential operators}. We refer to them as \textit{pseudo-differential operators}. Then, following \cite{Lindgren2011}, we will make a slight abuse of language and we will call 
SPDEs the class of stochastic equations considered in this work. This class of SPDEs is defined through operators of the form
\begin{equation}
\label{eq:DefLg}
\mathcal{L}_{g}( \cdot ) := \mathscr{F}^{-1}\left( g\mathscr{F}\left( \cdot \right)  \right),
\end{equation}
where $g: \mathbb{R}^{d} \to \mathbb{C}$ must be a sufficiently regular and Hermitian-symmetric function, that is it must satisfy $g(x) = \overline{g(-x)}$, where $\overline{a}$ is the complex conjugate of $a$. Under these conditions,  $\mathcal{L}_{g}$ is a real operator thanks to the properties of the Fourier transform. In this work, we will require that $g$ is continuous and bounded by a polynomial; see Section  \ref{Section:TheoreticalFramework} for a detailed exposition of all technical requirements. From now, every continuous, polynomially bounded and Hermitian-symmetric complex function $g$ defined over $\mathbb{R}^{d}$ will be called a \textit{symbol function} over $\mathbb{R}^{d}$, and it will be said to be the symbol function of the operator $\mathcal{L}_{g}$. Before presenting our main result, we first need to establish the relationship between the spectral measures of a stationary GeRF $U$ and its transform through  the operator $\mathcal{L}_g$ defined in (\ref{eq:DefLg}). The next proposition will be proven in Section \ref{Section:TheoreticalFramework}. 
\begin{Prop}
\label{PropositionLgUStationary}
Let $U$ be a real stationary GeRF over $\mathbb{R}^d$ with spectral measure $\mu_{U}$ and covariance distribution $\rho_{U}$, and let $g$ be a symbol function over $\mathbb{R}^{d}$. Then, $\mathcal{L}_{g}U$, where $\mathcal{L}_{g}$ is defined in  Eq. (\ref{eq:DefLg}), is a real stationary GeRF with spectral measure $\mu_{\mathcal{L}_{g}U} = |g|^{2}\mu_{U}$ and with covariance distribution $\rho_{\mathcal{L}_{g}U} = \mathcal{L}_{|g|^{2}}\rho_{U} $. 
\end{Prop}

\subsection{Statement of the main result}
\label{Section:StatementMainResult}
Let us consider a symbol function $g$ over $\mathbb{R}^{d}$ and a stationary GeRF $X$, which will be called from now on the \textit{source term}. A question that arises is to establish under which conditions on $g$ and $X$ the SPDE
\begin{equation}
\label{eq:LgU=X}
\mathcal{L}_{g}U = X,
\end{equation}
has a stationary solution, whether it is unique or not and, when solutions exist, whether we can characterize their covariance structures. Theorem \ref{TheoremExistAndUniqueStationary} provides a general answer to this question in the second order sense. That is, we shall only impose that the two sides of Eq. (\ref{eq:LgU=X}) have the same (generalized) covariance, which we write
\begin{equation}
\label{eq:LgU=X2moment}
\mathcal{L}_{g}U \stackrel{2nd \ o.}{=} X.
\end{equation}
This is not equivalent to require that $U$ solves (\ref{eq:LgU=X}) \textit{strictly}.  Under this more restrictive requirement, the evaluations of $\mathcal{L}_{g}U$ and $X$ over the same test functions (or at the same points in the case of random functions) are almost surely equal random variables. In the language of stochastic processes, this is equivalent to require that $\mathcal{L}_{g}U$ is a \textit{modification} of $X$. From a direct application of Proposition \ref{PropositionLgUStationary}, we get that a spectral measure $\mu_{U}$ of a potential stationary solution to (\ref{eq:LgU=X2moment}) must satisfy
\begin{equation}
\label{eq:g2muU=muX}
|g|^{2} \mu_{U} =  \mu_{X}.
\end{equation}
This kind of problem is called a \textit{division problem} in distribution theory. The existence of real stationary solutions to (\ref{eq:LgU=X2moment}) arises from the existence of solutions to (\ref{eq:g2muU=muX}) which are even and in $\mathcal{M}_{SG}^{+}(\mathbb{R}^{d})$. The explicit result is now formally presented in Theorem \ref{TheoremExistAndUniqueStationary}, which is our main result.
\begin{Theo}
\label{TheoremExistAndUniqueStationary}
Let $X$ be a real stationary GeRF over $\mathbb{R}^{d}$ with spectral measure $\mu_{X}$. Let $g : \mathbb{R}^{d} \to \mathbb{C}$ be a symbol function, and let $\mathcal{L}_{g}$ be an operator as defined in (\ref{eq:DefLg}) with symbol $g$. Then, there exists a real stationary GeRF solution to the equation (\ref{eq:LgU=X2moment}) if and only if there exists $N \in \mathbb{N}$ such that
\begin{equation}
\label{eq:ConditionExistence}
\int_{\mathbb{R}^{d}} \dfrac{d\mu_{X}(\xi)}{|g(\xi)|^{2}(1+|\xi|^{2})^{N}} < \infty.
\end{equation}
In this case, the measure 
\begin{equation}
\label{eq:MuU=g-2MuX}
d\mu_{U}(\xi) = |g(\xi)|^{-2}d\mu_{X}(\xi)
\end{equation}
is a spectral measure, and any real stationary GeRF with spectral measure $\mu_{U}$ solves (\ref{eq:LgU=X2moment}). Moreover, $\mu_{U}$ is the unique solution in $\mathcal{M}_{SG}^{+}(\mathbb{R}^{d})$ to (\ref{eq:g2muU=muX}) if and only if $|g| > 0 $.
\end{Theo}

\begin{Rem}
\label{Rem:MUuFinite}
When $N=0$, i.e. if $|g|^{-2}$ is integrable with respect to the measure $\mu_{X}$,  the measure $\mu_U$ is finite and  the solution $U$ is thus a  mean-square continuous  random function. 
This case was studied in \cite{Whittle1963}, where it is mentioned that solutions corresponding to non finite measures $\mu_U$ still make sense in some framework, the theory of which was at that time not completely available. Our work can be seen as one possible answer to this note.
\end{Rem} 

\begin{Rem}
\label{Rem:symbolInfBoundedPol}
A Sufficient Condition for Existence and Uniqueness (SCEU), regardless of the source term $X$, is to require that $|g|$ is inferiorly bounded by the inverse of a strictly positive polynomial. In this case the operator $\mathcal{L}_{g} $ is actually invertible:  $1/g$ is a symbol function and it is straightforward that $\mathcal{L}_{1/g}$ is the inverse operator of $\mathcal{L}_{g} $. This implies that Eq. (\ref{eq:LgU=X}) can be solved explicitly with $ U = \mathcal{L}_{1/g}X $. Then,  by Proposition \ref{PropositionLgUStationary}, $U$ is the unique stationary solution and its spectral measure is (\ref{eq:MuU=g-2MuX}). We shall henceforth refer to this condition as the SCEU on $g$.
\end{Rem}

\begin{Rem}
\label{Rem:gWithZerosHomogeneousProblem}
When the closed set $g^{-1}(\lbrace 0 \rbrace) = \lbrace \xi \in \mathbb{R}^{d} \ | \ g(\xi)=0 \rbrace$ is non-empty, the non-uniqueness is due to the existence of stationary solutions to the homogeneous problem
\begin{equation}
\label{eq:LgUh=0}
\mathcal{L}_{g}U_{H} = 0 .
\end{equation}
Indeed, for a  spectral measure $\mu_{U_{H}}$ over $\mathbb{R}^{d}$ supported on  $g^{-1}(\lbrace 0 \rbrace) $, its associated stationary random field satisfies strictly Eq.  (\ref{eq:LgUh=0}) since $\mu_{{\cal L}_g U_H} =  |g|^{2}\mu_H = 0$. Thus, if existence is provided, the addition of any stationary solution to (\ref{eq:LgU=X2moment}) with a non-trivial independent stationary solution to (\ref{eq:LgUh=0}) is also a stationary solution to (\ref{eq:LgU=X2moment}), which implies non-uniqueness. 
\end{Rem}

\section{Theoretical framework and proof of the main result}
\label{Section:TheoreticalFramework}

In order to prove Theorem \ref{TheoremExistAndUniqueStationary}, it is necessary to lay out some theoretical background, which uses Schwartz's theory of Distributions and its application to construct GeRFs. We assume that the reader is familiar with the Schwartz space of test functions over $\mathbb{R}^d$, denoted ${\cal S}(\mathbb{R}^{d})$, its dual space of tempered distributions, ${\cal S}'(\mathbb{R}^{d})$, the space of multiplicators of the Schwartz space $\mathcal{O}_{M}(\mathbb{R}^{d})$, and the definition and properties of the Fourier transform $\mathscr{F}$ \citep{Schwartz1966}. For sake of completeness, essential reminders on tempered distributions are provided in Appendix \ref{Section:TemperedDistribution}. We suggest \cite{Ito1954} and \citet[chapter 10]{Matheron1965} for a more complete introduction to GeRFs. Proposition \ref{PropositionLgUStationary} and Theorem \ref{TheoremExistAndUniqueStationary} will be proven here. This Section can be skipped in a first reading by readers more interested in the spatial and spatio-temporal models of random fields obtained within this framework. 

\subsection{Slow-growing measures and pseudo-differential operators}
\label{Section:SlowGrowMeasuresAndPseudoDiffOperators}

A complex Radon measure $\mu$ over  $\mathbb{R}^{d}$ is said to be a \textit{slow-growing measure} if there exists $N \in \mathbb{N}$ such that the measure $(1+|x|^{2})^{-N}|\mu|$ is a finite measure, where $|\mu|$ denotes the measure of total variation of $\mu$; see \citet[chapter 6]{Rudin1987},  or \citet[chapter 1.A]{Demengel2000}. We shall denote ${\cal M}_{SG}(\mathbb{R}^{d})$ the set of all slow-growing complex measures over $\mathbb{R}^{d}$. Obviously, ${\cal M}_{SG}^{+}(\mathbb{R}^{d}) \subset {\cal M}_{SG}(\mathbb{R}^{d})$. For a measure $\mu \in {\cal M}_{SG}(\mathbb{R}^{d})$,  the integral $ \langle \mu , \varphi \rangle := \int_{\mathbb{R}^{d}}\varphi(x)d\mu(x)$ is  well-defined for all $\varphi \in \mathcal{S}(\mathbb{R}^{d})$ and one can prove that it defines a tempered distribution, thus $\mathcal{M}_{SG}(\mathbb{R}^{d}) \subset \mathcal{S}'(\mathbb{R}^{d})$.

Let $g : \mathbb{R}^{d} \to \mathbb{C}$ be a polynomially bounded continuous function. The multiplication of $g$ by a slow-growing measure $\mu \in {\cal M}_{SG}(\mathbb{R}^{d})$, noted $g\mu$, is defined by  $(g\mu)(A):=\int_{A}g(x)d\mu(x)$ for every bounded Borel set $A \subset \mathbb{R}^{d}$. One can prove that $g\mu \in \mathcal{M}_{SG}(\mathbb{R}^{d})$, and thus it defines a tempered distribution through the expression
\begin{equation}
\langle g\mu , \varphi \rangle = \langle \mu , g\varphi \rangle = \int_{\mathbb{R}^{d}}\varphi(x)g(x)d\mu(x) \quad \forall \varphi \in \mathcal{S}(\mathbb{R}^{d}).
\end{equation}

As a consequence, pseudo-differential operators of the form $\mathcal{L}_{g} = \mathscr{F}^{-1}(g\mathscr{F}(\cdot ))$ as defined in Eq. (\ref{eq:DefLg}), with $g$ being a symbol function, are well-defined within our framework whenever the Fourier transform of the argument is a slow-growing measure. The domain of definition of $\mathcal{L}_{g}$ is thus the space of all tempered distributions such that its Fourier transform is a slow-growing measure:
\begin{equation}
\label{eq:DomainLg}
D(\mathcal{L}_{g}) = \lbrace T \in \mathcal{S}'(\mathbb{R}^{d}) \ | \ \mathscr{F}(T) \in \mathcal{M}_{SG}(\mathbb{R}^{d}) \rbrace.
\end{equation}
This class of operators includes for example linear combinations of differential operators which correspond to $g$ being an Hermitian-symmetric polynomial. Some fractional-differential operators are also included by taking $g$ to be a suitable continuous functions. A list of specific examples will be  worked out in Sections \ref{Section:IllustrationsSpatial} and \ref{Section:SpatioTemporalModels}.

\subsection{Generalized random fields}

A real $L^{2}-$tempered random distribution $Z$, referred to as real \textit{Generalized Random Field }(GeRF) from now on, is a real and continuous linear application from $\mathcal{S}(\mathbb{R}^{d})$ to $L^{2}(\Omega , \mathcal{F} , \mathbb{P})$, for some probability space $(\Omega , \mathcal{F} , \mathbb{P})$. 
We will  write  $\langle Z , \varphi \rangle := Z(\varphi) $ to emphasize that $Z$ acts as a continuous linear functional. All linear operators that are well-defined for tempered distributions can be used without restrictions on GeRFs, since they are defined through actions on test functions. In particular, differentiation and Fourier transforms are admissible operations on GeRFs (see Appendix \ref{Section:TemperedDistribution} for their definitions in the deterministic case). 

If $Z$ is a real GeRF, there exists a real mean distribution $ m_{Z} \in \mathcal{S}'(\mathbb{R}^{d})$ and a real covariance distribution $ C_{Z} \in \mathcal{S}'(\mathbb{R}^{d}\times\mathbb{R}^{d})$ satisfying $\mathbb{E}(\langle Z , \varphi \rangle ) = \langle m_{Z} , \varphi \rangle $ and $\mathbb{C}ov( \langle Z , \varphi \rangle , \langle Z , \phi \rangle ) = \langle C_{Z} , \varphi \otimes \overline{\phi} \rangle $ respectively for all $\varphi , \phi \in \mathcal{S}(\mathbb{R}^{d})$. Without loss of generality,  we will assume  $m_{Z}=0$. The covariance distribution must be a positive-definite kernel, i.e. it must satisfy $\langle C_{Z} , \varphi \otimes \overline{\varphi} \rangle \geq 0$ for all $\varphi \in \mathcal{S}(\mathbb{R}^{d})$, where $\otimes$ denotes the tensor product: $ ( \varphi \otimes \phi )(x,y) = \varphi(x)\phi(y)$ for $\varphi, \phi \in \mathcal{S}(\mathbb{R}^{d})$. The existence of this covariance distribution, which does not follow obviously from our assumptions, can be guaranteed by the Schwartz's Kernel Theorem applied to the space $\mathcal{S}'(\mathbb{R}^{d})$. See \citet[Theorem V.12]{ReedAndSimon1980} or \citet[Theorem 51.6 and its corollary]{Treves1967}.

A real GeRF $Z$ is \textit{second order stationary} (from now on, more simply, \textit{stationary}) if there exists a real and even distribution $\rho_{Z} \in \mathcal{S}'(\mathbb{R}^{d})$ such that $ \langle C_{Z} , \varphi \otimes \overline{\phi} \rangle = \langle \rho_{Z} , \varphi \ast \check{\overline{\phi}} \rangle $, where $\ast$ denotes the convolution product and $\ \check{ } \ $ denotes the reflection operator: $\check{\phi}(x) = \phi(-x)$. The distribution $\rho_{Z}$ is the \textit{stationary covariance distribution} of $Z$, or more simply the \textit{covariance distribution} if stationarity is clear from the context. It must be positive-definite, i.e., it must satisfy $\langle \rho_{Z} ,\varphi \ast \check{\overline{\varphi}} \rangle \geq 0$ for all $\varphi \in \mathcal{S}(\mathbb{R}^{d})$. A generalization of Bochner's Theorem, known as Bochner-Schwartz Theorem \citep[see for example][chapter 42]{Donoghue1969}, allows to conclude that any positive-definite and even distribution is the Fourier transform of a positive and even slow-growing measure. Thus, for $\rho_{Z}$ there exists a unique even measure $\mu_{Z} \in \mathcal{M}_{SG}^{+}(\mathbb{R}^{d})$ such that $\rho_{Z} = \mathscr{F}(\mu_{Z})$, known as the \textit{spectral measure} of $Z$. Since both $\mu_{Z}$ and $\rho_{Z}$ are even distributions, we will use extensively the following fact: $\rho_{Z} = \mathscr{F}(\mu_{Z}) = \mathscr{F}^{-1}(\mu_{Z})$. 

\subsection{Slow-growing orthogonal random measures}

A (not necessarily real) GeRF $Z$ is said to be a \textit{slow-growing random measure} if its covariance distribution $C_{Z}$ is a slow-growing measure, i.e. if $C_{Z} \in {\cal M}_{SG}(\mathbb{R}^d\times\mathbb{R}^{d})$. Similarly to slow-growing measures, slow-growing random measures can be multiplied by (deterministic) polynomially bounded continuous functions, thereby defining a new slow-growing random measure with covariance distribution in ${\cal M}_{SG}(\mathbb{R}^d\times\mathbb{R}^{d})$.   Specifically we have the following Proposition which is proven in Appendix \ref{Section:ProofPropositiongZ}:
\begin{Prop}
\label{Prop:MultiplicationgZ}
Let $Z$ be a slow-growing random measure with covariance $C_{Z} \in {\cal M}_{SG}(\mathbb{R}^d\times\mathbb{R}^{d})$ and let $g:\mathbb{R}^{d} \to \mathbb{C}$ be a polynomially bounded continuous function. Let us define the multiplication $gZ$ as a GeRF determined by $\langle gZ , \varphi \rangle = \langle Z , g\varphi \rangle $ for all $\varphi \in \mathcal{S}(\mathbb{R}^{d})$. Then the multiplication $gZ$ is well-defined as a GeRF and it is a slow-growing random measure with $C_{gZ}=  (g\otimes \overline{g} )C_{Z} \in {\cal  M }_{SG}(\mathbb{R}^d\times\mathbb{R}^{d})$. 
\end{Prop}
A particular class of slow-growing random measures is the class of \textit{slow-growing orthogonal random measures}, characterized by covariances of the form
\begin{equation}
\langle C_{Z} , \varphi \otimes \overline{\phi} \rangle = \int_{\mathbb{R}^{d}}\varphi(x)\overline{\phi}(x)d\nu_{Z}(x), \quad \varphi , \phi \in \mathcal{S}(\mathbb{R}^{d}),
\end{equation}
with $\nu_{Z} \in {\cal M}_{SG}^+(\mathbb{R}^{d})$. This form is obtained when the covariance measure $C_{Z}$ is supported over the hyperplane $\lbrace (x,y) \in \mathbb{R}^{d}\times\mathbb{R}^{d} \mid x = y \rbrace$, and it is easy to prove that it defines a positive-definite kernel. The measure $\nu_Z$ is called \textit{the weight} of $Z$. An important characteristic of this class is that these random measures take non-correlated values when evaluated over test functions that are orthogonal with respect to the weight measure, in particular when they have disjoint supports. From Proposition \ref{Prop:MultiplicationgZ} we get directly the following corollary.
\begin{Corollary}
\label{CorolLemmaMultiplicationgZ}
Let $Z$ be an orthogonal random measure with weight $\nu_{Z}$ and let $g$ be a continuous and polynomially bounded function. Then, $gZ$ is an orthogonal random measure with weight $|g|^{2}\nu_{Z}$.
\end{Corollary}
Slow-growing orthogonal random measures are in close connection with stationary GeRFs. A well known result, which is also easy to prove within the framework of GeRFs \citep[see][chapter 10]{Ito1954,Matheron1965}, is that the Fourier transform of a real stationary GeRF with spectral measure $\mu_{Z}$ is a Hermitian-symmetric complex slow-growing orthogonal random measure with weight $(2\pi)^{d/2}\mu_{Z}$. Grounded on this result, the Fourier transform of a stationary GeRF can be seen as a slow-growing measure. Operators of the form (\ref{eq:DefLg}), defined trough a symbol $g$,  can therefore be applied without restrictions. Having laid out these theoretical foundations, we are now able to prove Proposition \ref{PropositionLgUStationary}.

\subsection{Proof of the main result}

\subsubsection{Proof of Proposition \ref{PropositionLgUStationary}}

Let $g$ be a symbol function and let $\mathcal{L}_{g}$ its associated operator. Let $U$ be a real stationary GeRF with spectral measure $\mu_{U}$ and covariance distribution $\rho_{U}$. We know that $\mathscr{F}(U)$ is a Hermitian-symmetric complex slow-growing orthogonal random measure with weight $(2\pi)^{d/2}\mu_{U}$. Thus, by Corollary \ref{CorolLemmaMultiplicationgZ}, its multiplication by $g$ is well-defined and  is also a slow-growing orthogonal random measure with weight $(2\pi)^{d/2}|g|^{2}\mu_{U} \in \mathcal{M}_{SG}^{+}(\mathbb{R}^{d})$. Moreover, it is Hermitian-symmetric since $g$ is a symbol function. Hence, the inverse Fourier transform of $g\mathscr{F}(U)$, which is equal to $\mathcal{L}_{g}U$, is a real stationary GeRF with spectral measure $|g|^{2}\mu_{U}$. The expression of the covariance of $\mathcal{L}_{g}U$ is obtained immediately from $\rho_{\mathcal{L}_{g}U} = \mathscr{F}^{-1}\left(|g|^{2}\mu_{U}\right) = \mathscr{F}^{-1}\left(|g|^{2}\mathscr{F}(\rho_{U})\right) = \mathcal{L}_{|g|^{2}}\rho_{U} $.  $\blacksquare$

\subsubsection{Proof of Theorem \ref{TheoremExistAndUniqueStationary}}

Let $X$ be a real stationary GeRF over $\mathbb{R}^{d}$ with spectral measure $\mu_{X}$. Let $g$ be a symbol function over $\mathbb{R}^{d}$ and let $\mathcal{L}_{g}$ be its associated operator. We start by proving the existence criterion. Let us prove the necessity. Suppose there exists a real stationary GeRF, say $U$, satisfying (\ref{eq:LgU=X2moment}). By Proposition \ref{PropositionLgUStationary}, this implies that $|g|^{2}\mu_{U} = \mu_{X} $, and in particular we have that $\mu_{X}(  g^{-1}(\lbrace 0 \rbrace) ) = 0 $. As $\mu_{U} \in \mathcal{M}_{SG}^{+}(\mathbb{R}^{d})$, we can take $N \in \mathbb{N}$ such that $\int_{\mathbb{R}^{d}} (1+|\xi|^{2})^{-N}d\mu_{U}(\xi) < \infty $. We have that
\small
\begin{equation}
\int_{\mathbb{R}^{d}} \dfrac{d\mu_{X}(\xi)}{(1+|\xi|^{2})^{N}|g(\xi)|^{2}} = \int_{\lbrace g \neq 0 \rbrace} \dfrac{|g(\xi)|^{2}}{(1+|\xi|^{2})^{N}}\dfrac{d\mu_{U}(\xi)}{|g(\xi)|^{2}} = \int_{\lbrace g \neq 0 \rbrace} \dfrac{d\mu_{U}(\xi)}{(1+|\xi|^{2})^{N}} \leq  \int_{\mathbb{R}^{d}} \dfrac{d\mu_{U}(\xi)}{(1+|\xi|^{2})^{N}}  < \infty.
\end{equation}
\normalsize
Let us prove the sufficiency. The condition (\ref{eq:ConditionExistence}) implies in particular that the function $|g|^{-2}$ is locally integrable with respect to $\mu_{X}$. We can therefore define the Radon measure  $\mu_{U}(A) := \int_{A}|g(\xi)|^{-2}d\mu_{X}(\xi)$, for any bounded Borel set $A \subset \mathbb{R}^{d}$. By (\ref{eq:ConditionExistence}) and by the fact that both $\mu_{X}$ and $|g|^{2}$ are even, we see in addition that $\mu_{U} \in \mathcal{M}_{SG}^{+}(\mathbb{R}^{d})$ and it is even. Therefore $\mu_{U}$ is a spectral measure. Condition (\ref{eq:ConditionExistence}) also implies that $\mu_{X}(  g^{-1}(\lbrace 0 \rbrace) ) = 0 $. It is therefore straightforward that $|g|^{2}\mu_{U} = \mu_{X}$. Thus, any real stationary GeRF with spectral measure $\mu_{U}$ satisfies (\ref{eq:LgU=X2moment}).

Let us now prove the uniqueness criterion. For the necessity, suppose that $g$ does have zeros. Let us consider $\mu_{H}$, a slow-growing positive measure supported on the closed manifold $ g^{-1}(\lbrace 0 \rbrace) $. For instance, we can take any point $\xi_{0} \in \mathbb{R}^{d}$ such that $g(\xi_{0}) = 0 $ and use $\mu_{H} = \delta_{\xi_{0}} + \delta_{-\xi_{0}} $, which is a spectral measure. Hence, we have that $|g|^{2}\mu_{H} = 0 $. Thus, $\mu_{H}$ can be added to any solution $\mu_{U}$ of (\ref{eq:g2muU=muX}) and we will still get $|g|^{2}(\mu_{U} + \mu_{H}) = \mu_{X}$. We conclude that the solution is not unique. For the sufficiency, suppose $|g| > 0$ and that there are two different spectral measures $\mu_{1}$ and $\mu_{2}$ satisfying (\ref{eq:MuU=g-2MuX}). Then, the signed measure $\mu = \mu_{1} - \mu_{2}$ satisfies $|g|^{2}\mu = 0 $, and thus for any continuous function with compact support $\varphi$ we have $ \langle |g|^{2}\mu , \varphi \rangle = 0$. As $|g|$ is continuous and strictly positive, $|g|^{-2}\varphi$ is also continuous with compact support, and we can argue that for all $\varphi$ continuous with compact support,
\begin{equation}
\langle \mu , \varphi \rangle = \langle \mu , |g|^{2}|g|^{-2}\varphi \rangle = \langle |g|^{2}\mu , |g|^{-2}\varphi \rangle = 0.
\end{equation}
We conclude that $\mu = 0$ necessarily, and so $\mu_{1} = \mu_{2}$ and the solution is unique. $\blacksquare$

\section{White Noise as a fundamental case: a convolution theorem}
\label{Section:WhiteNoiseSourceTerm}

Let introduce the \textit{White Noise}, denoted $W$, defined as a real GeRF whose covariance distribution over $\mathbb{R}^{d}\times\mathbb{R}^{d}$ is
\begin{equation}
\label{eq:CovarianceWhiteNoise}
\langle C_{W} , \varphi \otimes \overline{\phi} \rangle = \int_{\mathbb{R}^{d}}\varphi(x)\overline{\phi}(x)dx, \quad \varphi , \phi \in \mathcal{S}(\mathbb{R}^{d}).
\end{equation}
$W$ is stationary with covariance distribution $\rho_{W} = \delta $, where $\delta \in {\cal S}'(\mathbb{R}^{d})$ is the Dirac measure at $0$. Its spectral measure is then proportional to the Lebesgue measure, $d\mu_{W}(x) = (2\pi)^{-\frac{d}{2}}dx$. $W$ is also a particular case of an orthogonal random measure whose weight is the Lebesgue measure. Since $\mu_{W}$ is not a finite measure, $W$ is not a random function and it can only be defined as a GeRF or as a random measure. We will see that SPDEs with a White Noise source term correspond to a  \textit{fundamental case} that can be used to obtain the covariance of solutions with more general source terms. Let us consider the SPDE
\begin{equation}
\label{eq:LgU = W2moment}
\mathcal{L}_{g}U \stackrel{2nd \ o.}{=} W.
\end{equation}
Theorem \ref{TheoremExistAndUniqueStationary} allows us to conclude that there are stationary solutions of (\ref{eq:LgU = W2moment}) if and only if the measure $|g|^{-2}(\xi)d\xi$ is in $\mathcal{M}_{SG}^{+}(\mathbb{R}^{d})$. We suppose this holds and we note $d\mu_{U}^{W}(\xi) = (2\pi)^{-\frac{d}{2}}|g|^{-2}(\xi)d\xi$, and $\rho_{U}^{W} = \mathscr{F}(\mu_{U}^{W})$. According to Proposition \ref{PropositionLgUStationary}, Eq. (\ref{eq:LgU = W2moment}) implies that
\begin{equation}
\label{eq:RhoUWGreenFunction}
\mathcal{L}_{|g|^{2}}\rho_{U}^{W} = \rho_{W} = \delta.
\end{equation}
Hence, the covariance $\rho_{U}^{W}$ can be seen as a \textit{Green's Function} of the operator $\mathcal{L}_{|g|^{2}}$. It turns out that in order to find the covariance of a solution to (\ref{eq:LgU=X2moment}) with an arbitrary source term $X$, we have to study the convolvability between $ \rho_{U}^{W}$ and $\rho_{X}$. If convolvability is satisfied, we get $\rho_{U} = \rho_{U}^{W} \ast \rho_{X}$. Theorem \ref{Theorem:CovarianceCovolutionWhiteNoiseSourceTerm} provides a sufficient criteria regarding the applicability of this procedure regardless of the source term $X$.

\begin{Theo}
\label{Theorem:CovarianceCovolutionWhiteNoiseSourceTerm}
Let $X$ be a real stationary GeRF over $\mathbb{R}^{d}$ with covariance distribution $\rho_{X}$. Let $g$ be a symbol function over $\mathbb{R}^{d}$ such that $\frac{1}{g}$ is smooth with polynomially bounded derivatives of all orders. Then, there exists a unique stationary solution to (\ref{eq:LgU=X2moment}) and its covariance distribution is given by
\begin{equation}
\label{eq:CovarianceUConvolutionXWhiteNoise}
\rho_{U} = \rho_{U}^{W} \ast \rho_{X},
\end{equation}
where $\rho_{U}^{W}$ is the covariance of the unique stationary solution to (\ref{eq:LgU = W2moment}).
\end{Theo}

\noindent \textbf{Proof: } Since $1/g$ is smooth and polynomially bounded, the SCEU holds, and there exists a unique stationary solution to equation (\ref{eq:LgU=X2moment}). The spectral measure of the solution is given by $\mu_{U} = |g|^{-2}\mu_{X} $,  $\mu_{X}$ being the spectral measure of $X$. The regularity and boundedness conditions for $1/g$ and its derivatives imply that both $\frac{1}{g}$ and $|g|^{-2}$ are in the space $\mathcal{O}_{M}(\mathbb{R}^{d})$ of multiplicators of the Schwartz space (see Appendix \ref{Section:TemperedDistribution}). Since the expression $|g|^{-2}\mu_{X}$ is the multiplication between $|g|^{-2} \in \mathcal{O}_{M}(\mathbb{R}^{d})$ and $\mu_{X} \in \mathcal{S}'(\mathbb{R}^{d})$,  the exchange formula for the Fourier transform can be applied. We thus obtain
\begin{equation}
\rho_{U} = \mathscr{F}(\mu_{U})= \mathscr{F}( |g|^{-2}\mu_{X} )  = \mathscr{F}( (2\pi)^{-\frac{d}{2}}|g|^{-2} ) \ast \mathscr{F}(\mu_{X})= \mathscr{F}(\mu_{U}^{W} ) \ast \mathscr{F}(\mu_{X}) = \rho_{U}^{W} \ast \rho_{X} . \quad  \blacksquare
\end{equation}
Although the condition on $g$ required in Theorem \ref{Theorem:CovarianceCovolutionWhiteNoiseSourceTerm} may seem restrictive, it turns out that it is satisfied by most models studied in the statistical literature on spatio-temporal random fields. For example, the Matérn model, Markov models and the Stein model satisfy these conditions, as it will be detailed in Sections \ref{Section:MaternModel}, \ref{Section:MarkovModel} and \ref{Section:SteinModel}. A more general analysis could be done by studying the convolvability between $\rho_{U}^{W}$ and $\rho_{X}$ in a classical framework as convolution between functions when applicable, or in the more general framework of the $\mathcal{S}'$--convolution; see e.g. \cite{DierolfAndVoigt1978}. 

Theorem \ref{Theorem:CovarianceCovolutionWhiteNoiseSourceTerm} shows that solutions of SPDEs with White Noise source term is the starting point of more general solutions, when the source term can be any stationary GeRF. In the next Sections, devoted to spatial and spatio-temporal models we shall always give special importance to the case of a White Noise source term.

\section{Application to known spatial models}
\label{Section:IllustrationsSpatial}

\subsection{Matérn Model}
\label{Section:MaternModel}
As a first example, we  start with a well-known and increasingly popular model, namely the Matérn model. The relationship between the Matérn Model and the SPDE over $\mathbb{R}^{d}$
\begin{equation}
\label{eq:SPDEMatern}
(\kappa^{2} - \Delta)^{\frac{\alpha}{2}}U = W,
\end{equation}
with $\kappa > 0 $, $\alpha \in \mathbb{R}$ and where $\Delta$ denotes the Laplace operator, has been established a long time ago \citep{Whittle1963} and recently revisited in \citet{Lindgren2011}. It can be easily be re-obtained from Theorem \ref{TheoremExistAndUniqueStationary}. The operator $(\kappa^{2} - \Delta)^{\frac{\alpha}{2}}$ is of the form (\ref{eq:DefLg}) with symbol function $g(\xi) = (\kappa^{2} + |\xi|^{2})^{\frac{\alpha}{2}}$, satisfying the SCEU defined in Remark \ref{Rem:symbolInfBoundedPol}. This allows us to conclude that there exists a unique stationary solution to (\ref{eq:SPDEMatern}), with spectral measure
\begin{equation}
\label{eq:SpectralMeasureMatern}
d\mu_{U}^{W}(\xi) = \dfrac{d\xi}{(2\pi)^{\frac{d}{2}}(\kappa^{2} + |\xi|^{2})^{\alpha}}.
\end{equation}
If $\alpha > \frac{d}{2}$, the measure (\ref{eq:SpectralMeasureMatern}) is finite, and thus its associated random field is a mean-square continuous random function, with stationary Matérn covariance function
\begin{equation}
\label{eq:CovarianceMatern}
\rho_{U}^{W}(h) =\dfrac{1}{(2\pi)^{\frac{d}{2}}2^{\alpha-1}\kappa^{2\alpha-d}\Gamma(\alpha)}(\kappa |h|)^{\alpha-\frac{d}{2}}K_{\alpha-\frac{d}{2}}(\kappa |h|),
\end{equation}
where $\Gamma$ is the Gamma function and $K_{\alpha-\frac{d}{2}}$ is the modified Bessel function of the second kind of order $\alpha-\frac{d}{2}$. When $\alpha \leq \frac{d}{2}$, we still obtain a unique stationary solution, but it is only defined in a distributional sense. We refer to this covariance as the \textit{generalized Matérn covariance }.

Since $g$ also satisfies the conditions in Theorem \ref{Theorem:CovarianceCovolutionWhiteNoiseSourceTerm}, we get that for any real stationary GeRF $X$, the SPDE
\begin{equation}
\label{eq:SPDEMaternX}
(\kappa^{2} - \Delta)^{\frac{\alpha}{2}}U = X
\end{equation}
has a unique stationary solution whose covariance is the convolution between  $\rho_{X}$ and the generalized Matérn covariance.

\subsection{Matérn model without range parameter}
\label{Section:MaternWithoutKappa}

The condition $\kappa>0$ in the Matérn SPDE defined in Eq. (\ref{eq:SPDEMatern}) can be relaxed. Setting $\kappa = 0 $, we obtain a fractional Laplacian operator $(-\Delta)^{\frac{\alpha}{2}}$ with symbol function $g(\xi) = |\xi|^{\alpha}$ for $\alpha > 0 $. Let us consider the SPDE
\begin{equation}
\label{eq:SPDEMaternWithoutKappa}
(-\Delta)^{\frac{\alpha}{2}}U = W,
\end{equation}
which corresponds to the limit case of a Matérn model as $\kappa \to 0$. In Theorem \ref{TheoremExistAndUniqueStationary}, the existence condition (\ref{eq:ConditionExistence}) requires that there exists $N \in \mathbb{N}$ such that $\int_{\mathbb{R}^{d}} (1+| \xi |^{2})^{-N}|\xi|^{-2\alpha}d\xi  < \infty$. Because of the singularity at the origin, this is only possible if $\alpha < d/2$. In this case, the spectral measure of a particular stationary solution to the equation (\ref{eq:SPDEMaternWithoutKappa}) is
\begin{equation}
\label{eq:SpectralMeasureMaternWithoutKappa}
d\mu_{U}(\xi) = \frac{1}{(2\pi)^{\frac{d}{2}}}\frac{d\xi}{|\xi|^{2\alpha}}.
\end{equation}
The associated covariance distribution is its Fourier transform, which is the locally integrable function (see \cite{Donoghue1969}, chapter 32)
\begin{equation}
\label{eq:CovarianceMaternWithoutKappa}
\rho_{U}(h) = \frac{1}{\pi^{\frac{d}{2}}}\frac{\Gamma(\frac{d}{2}-\alpha)}{\Gamma(\alpha)}\dfrac{1}{|h|^{d-2\alpha}}.
\end{equation}
Note that the function $\rho_{U}$ in (\ref{eq:CovarianceMaternWithoutKappa}) is not defined at $h = 0$.  It is not continuous, but it is still positive-definite in distributional sense. The associated random field must be interpreted as a GeRF and not as a mean-square continuous random function. This is an example of the kind of covariance structures we obtain when working with non-finite spectral measures. Such models have a long-range dependence behavior. They have been studied in \cite{Anh1998} and in \cite{GayAndHeyde1990}, in which Eq. (\ref{eq:SPDEMaternWithoutKappa}) is specified with a slightly different definition of the operator $(-\Delta)^{\frac{\alpha}{2}}$.

We remark that the symbol function $g(\xi) = |\xi|^{\alpha}$ has a zero at the origin. Hence, the SCEU does not hold. The stationary solution associated to the covariance (\ref{eq:CovarianceMaternWithoutKappa}) is not the unique possible solution. To describe all possible stationary solutions, we follow Remark \ref{Rem:gWithZerosHomogeneousProblem} and we consider spectral measures which are supported at the origin, i.e., which are proportionals to the Dirac measure $\mu_{U_{H}} = a\delta $, with $a > 0 $. The associated covariance distributions are then constant positive functions, and thus the associated GeRF are \textit{random constants}, that is, $U_{H}(x) = A$, for all $x \in \mathbb{R}^{d}$, with $A$ being a centered random variable with variance $(2\pi)^{-\frac{d}{2}}a< \infty$. In other words, \textit{the only stationary solutions to the homogeneous equation $(-\Delta)^{\frac{\alpha}{2}}U_{H} = 0 $ are random constants}. 

\subsection{Isotropic Markov models}
\label{Section:MarkovModel}

Let $p : \mathbb{R}^{+} \to \mathbb{R}^{+}$ be a strictly positive polynomial over $\mathbb{R}^{+}$. We consider the SPDE over $\mathbb{R}^{d}$
\begin{equation}
\label{eq:SPDEIsotropicMarkov}
p^\frac{1}{2}(-\Delta)U = W,
\end{equation}
where the operator $p^\frac{1}{2}(-\Delta)$ is of the form (\ref{eq:DefLg}) with  symbol function $g(\xi) = p^\frac{1}{2}(|\xi|^{2}) $. The SCEU holds, and thus the SPDE (\ref{eq:SPDEIsotropicMarkov}) has a unique stationary solution, and its spectral measure is of the form
\begin{equation}
\label{eq:SpectralMeasureIsotropicMarkov}
d\mu_{U}^{W}(\xi) = \dfrac{1}{(2\pi)^{\frac{d}{2}}} \dfrac{d\xi}{p(|\xi|^{2})}.
\end{equation}
This is a measure whose density is the inverse of a strictly positive and isotropic polynomial. Rozanov's Theorem \citep{Rozanov1977} allows us to conclude that this model is an isotropic stationary Markov Random Field (MRF). According to Rozanov's Theory, a MRF is, broadly speaking, a GeRF such that for every domain of $\mathbb{R}^{d}$,  evaluations of the random field on the interior of the domain are independent upon evaluations on the interior of the complement of the domain, conditionally to the behavior of the random field on a neighborhood of the boundary of the domain. By \textit{evaluations}, we mean the action of the GeRF over test functions whose supports are in the interior of the corresponding set. Rozanov's theorem states that every stationary MRF has a spectral measure whose density is the inverse of a strictly positive polynomial. Thus, in the case of isotropic models, MRFs satisfy equation (\ref{eq:SPDEIsotropicMarkov}). See \citet{Rozanov1982} for a complete theory of MRFs which also uses the theory of GeRFs.

Note that $g$ satisfies the conditions of Theorem \ref{Theorem:CovarianceCovolutionWhiteNoiseSourceTerm}. Hence, for any real stationary GeRF $X$ there exists a unique stationary solution to the SPDE 
\begin{equation}
p^\frac{1}{2}(-\Delta)U = X,
\end{equation}
whose covariance is the convolution between $\rho_{X}$ and the covariance of the MRF solution to Eq. (\ref{eq:SPDEIsotropicMarkov}).

\section{Application to spatio-temporal models}
\label{Section:SpatioTemporalModels}

We now present stationary spatio-temporal models which can be obtained and described within our framework. From now on, $d$ will always denote the spatial dimension, and we will explicitly write $\mathbb{R}^{d}\times\mathbb{R}$ referring to the space-time domain. We will denote $\mathscr{F}, \mathscr{F}_{S}$ and $\mathscr{F}_{T}$, respectively the spatio-temporal, spatial and temporal Fourier transforms. We will use the variables $(\xi,\omega)  \in \mathbb{R}^{d}\times\mathbb{R}$ for the frequency space-time domain (that is, after applying a spatio-temporal Fourier transform). When working with stationary covariance functions or distributions, the spatial separation vector will always be denoted by $h \in \mathbb{R}^{d}$ and the temporal lag by $u \in \mathbb{R}$. The function $g$ will always denote  a \textit{spatial} symbol function. Thus,  $\mathcal{L}_{g}$ denotes the operator $\mathscr{F}_{S}^{-1}(g\mathscr{F}_{S}(\cdot))$, which will be applied  to stationary GeRFs over $\mathbb{R}^{d}\times\mathbb{R}$. We denote $g_{R}$ and $g_{I}$ the real and imaginary parts of $g$ respectively. As $g$ is Hermitian-symmetric, $g_{R}$ is even and $g_{I}$ is odd.
 
 We recall the important concepts of \textit{separability} and \textit{symmetry} of a spatio-temporal stationary model. A stationary GeRF $Z$ over $\mathbb{R}^{d}\times\mathbb{R}$ is said to be \textit{separable} if its covariance $\rho_{Z} \in \mathcal{S}'(\mathbb{R}^{d}\times\mathbb{R}) $ can be expressed as the tensor product of a spatial covariance and a temporal covariance, $\rho_{Z} = \rho_{Z_{S}}\otimes \rho_{Z_{T}}$, with $\rho_{Z_{S}} \in \mathcal{S}'(\mathbb{R}^{d})$ and $\rho_{Z_{T}} \in \mathcal{S}'(\mathbb{R})$, obtaining $\rho_{Z}(h,u) = \rho_{Z_{S}}(h)\rho_{Z_{T}}(u)$ in the case with functional meaning. This is equivalent to require the spatio-temporal spectral measure to be the tensor product between a spatial spectral measure and a temporal spectral measure,  $d\mu_{Z}(\xi , \omega) = d\mu_{Z_{S}}(\xi)d\mu_{Z_{T}}(\omega)$. If $Z$ is separable, we write $Z = Z_{S}\otimes Z_{T}$, $Z_{S}$ and $Z_{T}$ \textit{representing} the corresponding spatial and temporal GeRFs with covariance $\rho_{Z_{S}}$ and $\rho_{Z_{T}}$ respectively. A stationary GeRF $Z$ over $\mathbb{R}^{d}\times\mathbb{R}$ is said to be \textit{symmetric} if its covariance satisfies $\rho_{Z}(h,u) = \rho_{Z}(h,-u) = \rho_{Z}(-h,u) = \rho_{Z}(-h,-u) $ in the case with functional meaning, with its corresponding generalization in the case of distributions using reflections with respect to the corresponding components. Symmetry is equivalent to have $d\mu_{Z}(\xi,\omega) = d\mu_{Z}(-\xi,\omega) = d\mu_{Z}(\xi,-\omega) = d\mu_{Z}(-\xi,-\omega)$. Since this holds if and only if the measure $\mu_{Z}$ is invariant under reflection with respect to the variable $\omega$, a non-symmetric model is obtained when the measure $\mu_{Z}$  depends on $\omega$ not only trough its absolute value $|\omega|$.

The spatial behavior of a stationary random function over $\mathbb{R}^{d}\times\mathbb{R}$, say $Z$, with covariance function $\rho_{Z}$ is studied by fixing the time component at any particular time $t \in \mathbb{R}$, obtaining the spatial random function $Z_{S} = (Z(x,t))_{x\in\mathbb{R}^{d}}$. Because of  time stationarity, $Z_{S}$ has  the same spatial covariance function for any chosen $t$, with $\rho_{Z_{S}}(h) = \rho_{Z}(h,0)$. We refer to $Z_{S}$ as a \textit{spatial trace} of $Z$, and to $\rho_{Z_{S}}$ as the \textit{spatial margin} of $\rho_{Z}$.    This can be generalized to any stationary GeRF $Z$ such that its spectral measure is temporally integrable, that it, it satisfies  $\mu_{Z}(A\times\mathbb{R}) < \infty $ for every bounded Borel set $A \subset \mathbb{R}^{d}$. In such a case, $Z$ is continuous in time, and the spectral measure of a spatial trace of $Z$ is $\mu_{Z_{S}}(A) = (2\pi)^{-\frac{1}{2}}\mu_{Z}(A \times \mathbb{R})$ for every bounded Borel set $A \subset \mathbb{R}^{d}$. The covariance of the spatial traces is simply $\rho_{Z_{S}} = \mathscr{F}_{S}(\mu_{Z_{S}})$, and it is said to be the spatial margin of the distribution $\rho_{Z}$. 

\subsection{Stein model}
\label{Section:SteinModel}
 
 In this example, we start from the spectral measure over  $\mathbb{R}^{d}\times\mathbb{R}$ proposed in \cite{Stein2005}
 
 \begin{equation}
\label{eq:SpectralMeasureStein}
d\mu_{U}(\xi,\omega) = \frac{1}{(2\pi)^{\frac{d+1}{2}} }\dfrac{d\xi d\omega}{(b(s^{2}+\omega^{2})^{\beta} + a(\kappa^{2}+|\xi|^{2})^{\alpha}  )^{\nu}},
\end{equation}
with $a,b > 0 $, $s^{2} +\kappa^{2} > 0 $, and $\alpha , \beta , \nu \in \mathbb{R}$ and we derive a corresponding SPDE. The measure (\ref{eq:SpectralMeasureStein}) is always a well-defined spectral measure since its density is the inverse of a positive and polynomially bounded continuous function. We consider the spatio-temporal symbol function
\begin{equation}
\label{eq:SymbolStein}
(\xi,\omega) \mapsto (b(s^{2}+\omega^{2})^{\beta} + a(\kappa^{2}+|\xi|^{2})^{\alpha}  )^{\nu/2},
\end{equation}
which satisfies the SCEU. Using a spatio-temporal White Noise, $W$ ($d\mu_{W}(\xi,\omega) = (2\pi)^{-(d+1)/2}d\xi d\omega$), a corresponding SPDE for the Stein model is then
\begin{equation}
\label{eq:SPDEStein}
\left(b\left(s^{2}-\dfrac{\partial^{2}}{\partial t^{2}}\right)^{\beta} + a\left(\kappa^{2} - \Delta \right)^{\alpha} \right)^{\nu/2} U = W.
\end{equation}
As a consequence of Theorem \ref{TheoremExistAndUniqueStationary}, there exists  a unique stationary solution to ($\ref{eq:SPDEStein}) $ and its spectral measure is (\ref{eq:SpectralMeasureStein}). When $\alpha, \beta$ and $\nu$ are positive, and if $1/(\beta\nu) + d/(\alpha\nu) < 2 $ holds,  \citet{Stein2005} shows that the measure (\ref{eq:SpectralMeasureStein}) is finite and that its associated random field is a mean-square continuous random function. In other cases, the model is not a function but a distribution, and we refer to this model as the \textit{generalized Stein model.} The interesting property of the Stein model is that, without being a separable model,  the spatial and temporal smoothness of the paths of the random function can be controlled separately  thanks to the parameters $\alpha$ and $\beta$. Except for some particular values of the parameters, there is no closed-form expression for the covariance.

When $\kappa, s, a , b > 0 $ and $\alpha, \beta, \nu$ are not null, the symbol function (\ref{eq:SymbolStein}) satisfies conditions in Theorem \ref{Theorem:CovarianceCovolutionWhiteNoiseSourceTerm}. Hence, for any stationary GeRF $X$, the SPDE
\begin{equation}
\left(b\left(s^{2}-\frac{\partial^{2}}{\partial t^{2}}\right)^{\beta} + a\left(\kappa^{2} - \Delta \right)^{\alpha} \right)^{\nu/2} U = X
\end{equation}
has a unique stationary solution whose covariance is the convolution between $\rho_{X}$ and the covariance of the generalized Stein model.

\subsection{Models derived from evolution equations}
\label{Section:EvolutionModels}

In this section we study models associated to the following class of SPDEs over $\mathbb{R}^{d}\times\mathbb{R}$ which we call \textit{evolution equations}:
\begin{equation}
\label{eq:SPDEGeneralEvolutionModel}
\dfrac{\partial^{\beta}U}{\partial t^{\beta}} + \mathcal{L}_{g}U = X,
\end{equation}
where $X$ is a stationary spatio-temporal GeRF, and $\beta > 0 $. For this class of SPDEs, we study in detail several examples of physical and statistical interest. They involve for example advection, Langevin-type equation, heat diffusion and wave propagation phenomena.

First of all, for $\beta > 0$, we define the operator  $ \frac{\partial^{\beta}}{\partial t^{\beta}}$
\begin{equation}
\label{eq:DefFractionalDifferentialOperator}
\dfrac{\partial^{\beta}}{\partial t^{\beta}} := \mathscr{F}_{T}^{-1}( (i\omega)^{\beta}\mathscr{F}_{T}(\cdot) ), 
\end{equation}
where we have used the symbol function over $\mathbb{R}$
\begin{equation}
\label{eq:SymbolFractionalDifferentialOperator}
\omega \mapsto (i\omega)^{\beta} := |\omega|^{\beta}e^{i \sgn(\omega)\beta \pi/2}.
\end{equation}
The function (\ref{eq:SymbolFractionalDifferentialOperator}) is continuous, Hermitian-symmetric and bounded by a polynomial for every $\beta > 0 $, so it is a symbol function. Similar definitions of a fractional differential operator can be found in \cite{Maniardi2001}. We call a \textit{fractional order evolution model}  a spatio-temporal stationary solution of the SPDE (\ref{eq:SPDEGeneralEvolutionModel}) with $\beta \notin \mathbb{N}$. For $\beta \in \mathbb{N}$, (\ref{eq:DefFractionalDifferentialOperator}) coincides with a classical differential operator. The corresponding stationary solutions are called $n-$\textit{th order evolution model}.

The spatio-temporal symbol function of the operator involved in (\ref{eq:SPDEGeneralEvolutionModel}) is
\begin{equation}
\label{eq:SymbolFunctionGenOEM}
(\xi,\omega) \mapsto (i\omega)^{\beta} + g(\xi) = |\omega|^{\beta}\cos\left( \frac{\beta\pi}{2} \right) + g_{R}(\xi) + i\left(\sgn(\omega)|\omega|^{\beta}\sin \left(\frac{\beta\pi}{2} \right) + g_{I}(\xi) \right).
\end{equation}
Theorem \ref{TheoremExistAndUniqueStationary} allows us to conclude that there exists stationary solutions to (\ref{eq:SPDEGeneralEvolutionModel}) if the measure
\begin{equation}
\label{eq:SpectralMeasureGenOEM}
d\mu_{U}(\xi,\omega) = \frac{d\mu_{X}(\xi,\omega)}{|(i\omega)^{\beta} + g(\xi)|^{2}},
\end{equation}
is a well-defined spectral measure. We will follow Remark \ref{Rem:symbolInfBoundedPol} and look for conditions on $g$ such that (\ref{eq:SymbolFunctionGenOEM}) satisfies the SCEU and thus to have a unique stationary solution regardless of the source term $X$. The next proposition, proven in Appendix \ref{Section:ProofPropositionSCEUGenOEM}, allows us to identify the cases where the SCEU holds regardless of the imaginary part $g_{I}$.
\begin{Prop}
\label{Prop:SCEUGenOEM}
A necessary and sufficient condition for the function  (\ref{eq:SymbolFunctionGenOEM}) to satisfy the SCEU for every arbitrary $g_{I}$ function is that $g_{R}$ satisfies the SCEU and $g_R \cos(\frac{\beta\pi}{2}) \geq 0$.
\end{Prop}
We suppose that conditions over $g_{R}$ in  Proposition \ref{Prop:SCEUGenOEM} hold. Let us study the properties of this kind of models. We restrain ourselves to the cases where $X$ is a separable model $X = X_{S}\otimes X_{T}$. The spectral measure of the unique stationary solution to (\ref{eq:SPDEGeneralEvolutionModel}) is then
\begin{equation}
\label{eq:SpectralMeasureGenOEMXSeparable}
d\mu_{U}(\xi,\omega) = \frac{d\mu_{X_{S}}(\xi)d\mu_{X_{T}}(\omega)}{|\omega|^{2\beta} + 2 | \omega |^{\beta}\left( g_{R}(\xi)\cos\left( \frac{\beta\pi}{2} \right)  + \sgn(\omega)g_{I}(\xi)\sin\left( \frac{\beta \pi}{2}  \right)   \right) + | g(\xi) |^{2} }.
\end{equation}
A separable model is obtained when $g_{I} = 0$ and $g_{R}$ is a constant function. Otherwise, the model is not separable. The function $ \sgn $ in (\ref{eq:SpectralMeasureGenOEMXSeparable}) allows to identify the cases where the spectral measure does not depend on the argument $\omega$ only trough $|\omega|$ and thus the symmetry of the model can be controlled. A symmetric model is then obtained  when $\beta$ is an even integer or when the function $g_{I}$ is null. A non-symmetric model is obtained otherwise. In this case the non-symmetry can be parametrized by controlling the function $g_{I}$. The mean-square temporal regularity of the associated random field depends on the parameter $\beta$, as it can be seen by analyzing the temporal-integrability of the measure $\mu_{U}$. Thus, this model allows a practical control of the separability, symmetry and regularity conditions.

The behavior of a spatial trace of this model can be obtained if the measure $\mu_{U}$ is temporally-integrable, that is, if 
$ \int_{\mathbb{R}}  |(i\omega)^{\beta} + g(\xi)|^{-2}d\mu_{X_{T}}(\omega) < \infty $. Let us restrict ourselves to the case where $X$ is White Noise in time, $X = X_{S}\otimes W_{T}$, i.e. $d\mu_{X}(\xi , \omega) = d\mu_{X_{S}}(\xi)d\mu_{W_{T}}(\omega) = d\mu_{X_{S}}(\xi)(2\pi)^{-\frac{1}{2}}d\omega $. In that case, the measure $\mu_{U}$ is temporally-integrable when $\beta > \frac{1}{2}$. The spectral measure of the spatial traces can be obtained by calculating the corresponding integral. We set $g_{I}=0$. The general case $g_{I} \neq0 $ is much more technical. For the sake of a clear exposition, it is left aside in this general presentation. The spectral measure of a spatial trace $U_{S}$ is then
\begin{align}
d\mu_{U_{S}}(\xi) &  = \frac{1}{2\pi}\int_{\mathbb{R}}\frac{d\omega}{|\omega|^{2\beta} + 2|\omega|^{\beta}g_{R}(\xi)\cos\left(\frac{\pi}{2}\beta \right) + g_{R}^{2}(\xi) } d\mu_{X_{S}}(\xi)  \\
   & = \frac{|g_{R}(\xi)|^{\frac{1}{\beta}-2} }{\pi \beta} \underbrace{\int_{0}^{\infty} \frac{\theta^{\frac{1}{\beta} - 1}}{\theta^{2} + 2\theta\sgn(g_{R})\cos\left( \frac{\pi}{2}\beta \right) + 1 }d\theta}_{= I_{\beta}}d\mu_{X_{S}}(\xi),  \label{eq:ComputationSpatialTraceGenOEM}
\end{align}
where we have used the parity of the function with respect to $\omega$ and then used the change of variable $\omega = ( |g_{R}(\xi)|\theta )^{\frac{1}{\beta}} $. The integral $I_{\beta}$ does not depend on $\xi$ since $g_{R}$ does not change in sign. This integral can be computed \citep[see for instance][3.252.12]{Gradshteyn1994}. In particular, $I_{1} = I_{2} = \pi/2$. Then, the spatial traces of the solution satisfy the spatial SPDE
\begin{equation}
\label{eq:SPDESpatialTraceGenOEM}
\sqrt{\frac{\pi\beta}{I_{\beta}}}\mathcal{L}_{ |g_{R}|^{1 - \frac{1}{2\beta}} }U_{S} \stackrel{2nd \ o.}{=} X_{S}.
\end{equation}
This model has a continuous point-wise meaning when the function $|g_{R}|^{\frac{1}{\beta} - 2} $ is integrable with respect to the measure $\mu_{X_{S}}$, case in which the measure $\mu_{U}$ is a finite measure.

Theorem \ref{Theorem:CovarianceCovolutionWhiteNoiseSourceTerm} cannot be applied when  $\beta \notin \mathbb{N}$ since the symbol function (\ref{eq:SymbolFunctionGenOEM}) is not smooth. The case $\beta \in \mathbb{N}$ can be worked out supposing some regularity conditions on $g$. We present the corresponding analysis for the cases $\beta \in \lbrace 1 , 2 \rbrace$.

\medskip

A \textit{first order evolution model} is a stationary solution of  Eq. (\ref{eq:SPDEGeneralEvolutionModel}) when $\beta = 1$. Let us set $X = W$, the spatio-temporal White Noise. The spectral measure is then

\begin{equation}
\label{eq:SpectralMeasureFOEMWhiteNoise}
d\mu_{U}^{W}(\xi ,\omega) = \frac{1}{(2\pi)^{\frac{d+1}{2}}}\frac{d\xi d\omega}{(\omega + g_{I}(\xi))^{2} + g_{R}^{2}(\xi)}.
\end{equation}
From this we obtain that its covariance is of the form
\begin{equation}
\label{eq:CovarianceFOEMWhiteNoise}
\rho_{U}^{W}(h,u) = \mathscr{F}_{S}\left(  \xi \mapsto \frac{1}{(2\pi)^{\frac{d}{2}}}\dfrac{e^{i u g_{I}(\xi) - |u||g_{R}(\xi)|} }{ 2|g_{R}(\xi)| } \right)(h).
\end{equation}
This model can then be seen as a \textit{mixture of (complex) exponentials}. For ease of reading, we have used a functional notation for the variables $(h,u)$ in (\ref{eq:CovarianceFOEMWhiteNoise}), but $\rho_{U}^{W}$ is not necessarily a function.  Generally, it is  a distribution. A continuous function is obtained when $|g_{R}|^{-1}$ is an integrable function. The spatial margin of $\rho_{U}^{W}$ is obtained by setting $u=0$ in $\rho_{U}^{W}$. It does not depend on $g_{I}$. Thus, Eq. (\ref{eq:SPDESpatialTraceGenOEM}) can be used to describe the spatial behavior of the model for the case $X_{S} = W_{S}$, including the cases where $g_{I} \neq 0$. It is immediate that it can also be applied to an arbitrary $X_{S}$,  in which case the term $ d\xi$ in (\ref{eq:SpectralMeasureFOEMWhiteNoise}) is replaced by $(2\pi)^{d/2}d\mu_{X_{S}}(\xi)$. The spatial structure is thus completely described by $g_{R}$, while the spatio-temporal non-symmetry is described by $g_{I}$.  Theorem \ref{Theorem:CovarianceCovolutionWhiteNoiseSourceTerm} can be applied if $g_{R}$, $g_{I}$ and $1/g_{R}$ are in $\mathcal{O}_{M}(\mathbb{R}^{d}) $, since in this case the reciprocal of the spatio-temporal symbol function (\ref{eq:SymbolFunctionGenOEM}) is in $\mathcal{O}_{M}(\mathbb{R}^{d}\times\mathbb{R})$. We obtain in that case that the covariance of the solution with an arbitrary source term $X$ is the convolution $\rho_{U}^{W}\ast \rho_{X}$.

\medskip

A \textit{second order evolution model} is a stationary solution of  Eq. (\ref{eq:SPDEGeneralEvolutionModel}) when $\beta = 2$. Consider again $X = W$. The spectral measure is then
\begin{equation}
\label{eq:SpectralMeasureSOEMWhiteNoise}
d\mu_{U}^{W}(\xi , \omega) = \dfrac{1}{(2\pi)^{\frac{d+1}{2}}} \dfrac{d\xi d\omega}{(\omega^{2}-g_{R}(\xi))^{2} + g_{I}^{2}(\xi)  },
 \end{equation}
 and the covariance distribution $\rho_{U}^{W}$ is the Fourier transform of $\mu_{U}^{W}$. To simplify the notation, consider the non-null complex spatial function
\begin{equation}
\label{eq:gammaFunction}
\gamma(\xi) = \sqrt{\frac{|g(\xi)| + g_{R}(\xi)}{2}} + i\sqrt{\frac{ | g(\xi) | - g_{R}(\xi)}{2}}.
\end{equation}
Let us denote $\gamma_{R}$ and $\gamma_{I}$ the real and imaginary parts of $\gamma$ respectively. The covariance $\rho_{U}^{W}$ is then
\small
\begin{equation}
\label{eq:CovarianceSOEMWhiteNoise}
\rho_{U}^{W}(h,u) = \mathscr{F}_{S}\left(\xi \mapsto \frac{e^{-(|\gamma_{I}(\xi)| + i \gamma_{R}(\xi) )|u|  }}{(2\pi)^{\frac{d}{2}}8|\gamma_{I}(\xi)|^{2}} \left[  \frac{1}{|\gamma_{I}(\xi)| + i  \gamma_{R}(\xi) } +  \frac{e^{i2\gamma_{R}(\xi)|u|}}{|\gamma_{I}(\xi)| - i  \gamma_{R}(\xi) } + \frac{e^{i2\gamma_{R}(\xi)|u|} - 1}{i \gamma_{R}(\xi)}   \right]   \right)(h)
\end{equation}
\normalsize
where  $(e^{i2\gamma_{R}(\xi)|u|} - 1)/i \gamma_{R}(\xi) $ is set to  $2|u|$ when $\gamma_{R}(\xi) = 0 $, which corresponds to $g_{I}(\xi) = 0 $.  This covariance distribution is a continuous function if the function $|\gamma_{I}|^{-1}| \gamma |^{-2}$ is integrable over $\mathbb{R}^{d}$, which is equivalent to require that the function $|g|^{-1}(|g|-g_{R})^{-\frac{1}{2}} $ is integrable over  $\mathbb{R}^{d}$. Contrarily to the case of first order evolution models, this model is always symmetric and the structure of the spatial traces depends on both $g_{R}$ and $g_{I}$, as it can be seen by evaluating (\ref{eq:CovarianceSOEMWhiteNoise}) at $u=0$. Thus, Eq. (\ref{eq:SPDESpatialTraceGenOEM}) does not hold for $g_{I} \neq 0$. The spectral measure of the spatial traces is
\begin{equation}
\label{eq:SpectralMeasureSpatialTraceSOEMWhiteNoise}
d\mu_{U_{S}}^{W}(\xi) = \frac{d\xi}{(2\pi)^{\frac{d}{2}}2\sqrt{2} |g(\xi)|\sqrt{|g(\xi)|-g_{R}(\xi))} },
\end{equation}
from which we obtain that a spatial trace $U_{S}$ satisfies the spatial SPDE
\begin{equation}
\label{eq:SPDESpatialTraceSOEMWhiteNoise}
\sqrt{2\sqrt{2}}\mathcal{L}_{\sqrt{|g|\sqrt{|g|-g_{R}}}}U_{S} \stackrel{2nd \ o.}{=} W_{S},
\end{equation}
where $W_{S}$ is a spatial White Noise. An analogue expression is obtained in the case $X = X_{S}\otimes W_{T}$, by replacing $W_{S}$ by $X_{S}$ in (\ref{eq:SPDESpatialTraceSOEMWhiteNoise}) and $d\xi$ by $(2\pi)^{d/2}d\mu_{X_{S}}(\xi)$ in (\ref{eq:SpectralMeasureSOEMWhiteNoise}) and (\ref{eq:SpectralMeasureSpatialTraceSOEMWhiteNoise}). When $X$ is a general spatio-temporal stationary GeRF, a sufficient condition to apply Theorem \ref{Theorem:CovarianceCovolutionWhiteNoiseSourceTerm} is that $g_{R}$, $1/g_{R} $ and $g_{I}$ are in the space $\mathcal{O}_{M}(\mathbb{R}^{d}) $. In this case, the only stationary solution to the SPDE (\ref{eq:SPDEGeneralEvolutionModel}) with $\beta = 2$ has a covariance of the form $\rho_{U} = \rho_{U}^{W} \ast \rho_{X} $, where $\rho_{U}^{W}$ is given by  (\ref{eq:CovarianceSOEMWhiteNoise}).

We now present some particular models inspired by physical and statistical literature. In some cases Proposition \ref{Prop:SCEUGenOEM} can be applied. In other cases, there is no uniqueness and sometimes not even existence of stationary solutions.

\begin{Example}{\textbf{Evolving Matérn model.}}
\label{Ex:FOEMMatern}
\end{Example}
\noindent In the most general term, we call an \textit{evolving Matérn model} a stationary solution of the evolution equation Eq. (\ref{eq:SPDEGeneralEvolutionModel}) such that its spatial traces follow a Matérn model. Evolving Matérn models can be obtained by adequately controlling $g$, the structure of $X$ or both. Here, we restrict ourselves to  models that are solutions to equations of the form
\begin{equation}
\label{eq:SPDEEvolutionMatern}
\dfrac{\partial^{\beta}U}{\partial t^{\beta}} + s_{\beta} a (\kappa^{2} - \Delta)^{\frac{\alpha}{2}}U = W,
\end{equation}
where $W$ is as usual the spatio-temporal White Noise, $\kappa^{2} , a > 0, \alpha \in \mathbb{R}$, and $s_{\beta}$  is a parameter that takes the value $1$ or $-1$  depending conveniently on $\beta$ in order to obtain the conditions in Proposition \ref{Prop:SCEUGenOEM} for $g(\xi) =  s_{\beta} a (\kappa^{2} + | \xi |^{2} )^{\frac{\alpha}{2}}$.  There is then a unique stationary solution to (\ref{eq:SPDEEvolutionMatern}). Its spectral measure is
\begin{equation}
\label{eq:SpectralMeasureEvolutionMatern}
d\mu_{U}(\xi,\omega) = \frac{1}{(2\pi)^{\frac{d+1}{2}}}\frac{d\xi d\omega}{|\omega|^{2\beta} + 2 | \omega |^{\beta} a (\kappa^{2} + | \xi |^{2} )^{\frac{\alpha}{2}}|\cos\left( \frac{\beta\pi}{2} \right)|  + a^{2} (\kappa^{2} + | \xi |^{2} )^{\alpha} }.
\end{equation}
Following Eq. (\ref{eq:SPDESpatialTraceGenOEM}), when $\beta > \frac{1}{2}$ the spatial traces of this model follow the spatial SPDE
\begin{equation}
\label{eq:SPDESpatialTraceEvolutionMatern}
\sqrt{\frac{\pi\beta}{I_{\beta}}}a^{1-\frac{1}{2\beta}} (\kappa^{2}-\Delta)^{\frac{\alpha}{2}\left(1 - \frac{1}{2\beta}\right) }U_{S} \stackrel{2nd \ o.}{=} W_{S},
\end{equation}
where $W_{S}$ is a spatial White Noise. Direct identification between (\ref{eq:SPDESpatialTraceEvolutionMatern}) and (\ref{eq:SPDEMatern}) indicates that the spatial covariance  is thus a Matérn covariance. It has a functional meaning when $ \alpha \left( 1 -1/(2\beta) \right) > d/2$. Explicit expressions of the covariances can be  obtained using  the Fourier transform of radial functions \citep[chapter 41]{Donoghue1969}.  

\medskip

In particular, for $\beta = 1$, we get 
\begin{equation}
\label{eq:CovarianceFOEMMaternIntegral}
\rho_{U}(h,u) = \dfrac{1}{(2\pi)^{\frac{d}{2}}|h|^{\frac{d-2}{2}}} \int_{0}^{\infty} J_{\frac{d-2}{2}}(|h|r)\dfrac{e^{- a(\kappa^{2}+r^{2})^{\frac{\alpha}{2}}|u|  } }{ 2a(\kappa^{2}+r^{2})^{\frac{\alpha}{2}  } }r^{\frac{d}{2}}dr,
\end{equation}
where $J_{b}$ denotes the Bessel function of the first kind of order $b$. This model has also been proposed in \citet{JonesAndZhang1997}, in which an approach similar to our framework was followed for first order evolution equations. This is a symmetric non-separable model which can be identified as a mixture of  a $J-$Bessel model in space with an exponential model in time. Notice that in this case we could add a non-null imaginary part $g_{I}$ to the symbol function without changing the spatial behavior, thereby generating non-symmetric evolving Matérn models. However, in this case the expression (\ref{eq:CovarianceFOEMMaternIntegral}) no longer applies.

\medskip

For $\beta = 2$, one gets
\begin{equation}
\label{eq:CovarianceSOEMMaternIntegral}
\rho_{U}(h,u) = \dfrac{1}{(2\pi)^{\frac{d}{2}}|h|^{\frac{d-2}{2}}} \int_{0}^{\infty} J_{\frac{d-2}{2}}(|h|r)\dfrac{e^{-\sqrt{a} (\kappa^{2}+r^{2})^{\frac{\alpha}{4}}|u|  }(1+\sqrt{a}(\kappa^{2}+r^{2})^{\frac{\alpha}{4}}|u|) }{ 4a\sqrt{a}(\kappa^{2}+r^{2})^{\frac{3\alpha}{4}  } }r^{\frac{d}{2}}dr.
\end{equation}
This covariance is a mixture of $J-$Bessel model in space and a Matérn model in time since the spectral measure (\ref{eq:SpectralMeasureSOEMWhiteNoise}) has the form of a Matérn spectral measure in $\omega$. Notice that this mixture property does not hold for $\beta \not\in \{1,2\}$. 

Notice that both $g_{R}$ and $1/g_{R}$ are in $\mathcal{O}_{M}(\mathbb{R}^{d})$. Thus, for $\beta \in \{1,2\}$ Theorem \ref{Theorem:CovarianceCovolutionWhiteNoiseSourceTerm} can be applied. In these cases, the covariance of the solution to an equation of the form (\ref{eq:SPDEEvolutionMatern}) with an arbitrary source term $X$ is  the convolution between (\ref{eq:CovarianceFOEMMaternIntegral}) for $\beta = 1$ (respectively (\ref{eq:CovarianceSOEMMaternIntegral}) for $\beta = 2 $) and  $\rho_{X}$.

We finally remark that for the cases $\beta \in \mathbb{N}$ these models are Stein models. See the correspondences between Eq. (\ref{eq:SpectralMeasureStein}) and Eq. (\ref{eq:SpectralMeasureEvolutionMatern}) in those cases, considering the temporal scale parameter $s=0$. 

\begin{Example}{\textbf{Advection-diffusion equation.}}
\label{Ex:AdvectionDiffusion}
\end{Example}

\citet{Sigrist2015} propose estimation methods and simulation algorithms for the unique stationary solution of the SPDE over $\mathbb{R}^{d}\times\mathbb{R}$:
\begin{equation}
\label{eq:SPDEAdvectionDiffusion}
\dfrac{\partial U}{\partial t} + \kappa^{2}U + v^{T}\nabla U -  \div( \Sigma \nabla U ) = X_{S}\otimes W_{T},
\end{equation} 
where $\kappa > 0 $ is a damping parameter, $v \in \mathbb{R}^{d}$ is a velocity and $\Sigma $ is a symmetric positive-definite matrix controlling non-isotropic diffusion. $W_{T}$ is a temporal White Noise and $X_{S}$ is a stationary spatial random field. This equation, known as the \textit{advection-diffusion equation}, is a particular first order evolution model. Its spatial symbol function is $g(\xi) = \kappa^{2} + \xi^{T}\Sigma\xi + iv^{T}\xi $, for which conditions in Proposition \ref{Prop:SCEUGenOEM} are satisfied. Without advection ($v=0$), this equation was studied in \cite{Whittle1963} in a non-generalized framework. \citet{Sigrist2015} consider a Matérn Model for $X_{S}$, with smoothness parameter equals to $1$, corresponding to  $\alpha = 2 $ in (\ref{eq:CovarianceMatern}) when $d=2$. The spatial behavior of this model is described by the SPDE (\ref{eq:SPDESpatialTraceGenOEM}) for $\beta = 1$. 

\begin{Example}{\textbf{A Langevin equation.}}
\label{Ex:LangevinEquation}
\end{Example}

Using linear response theory, \citet{Hristopulos2016} propose stationary random fields which are solutions to the Langevin equation 
\begin{equation}
\label{eq:SPDELangevin}
\dfrac{\partial U}{\partial t} + \frac{D}{2k^{d}\eta_{0}}\left( 1 - \eta_{1}k^{2}\Delta + \nu k^{4} \Delta^{2} \right) U = W,
\end{equation}
with $D, k, \eta_{0} > 0 $, $\eta_{1} ,\nu \geq 0 $. Let $C = D/(2k^d\eta_0)$. For this first order evolution model, the spatial symbol function is $ g(\xi) = C\left( 1 + \eta_{1}k^{2}|\xi|^{2} + \nu k^{2}|\xi|^{4} \right) $, which satisfies conditions of Proposition \ref{Prop:SCEUGenOEM}. Hence, (\ref{eq:SPDELangevin}) has a unique stationary solution, whose spectral measure can be obtained using the general expression of first order evolution model in (\ref{eq:SpectralMeasureFOEMWhiteNoise}).  \citet{Hristopulos2016} provide expressions of the related covariance structures, which are functions for $d \leq 3$, and which can be obtained through formulas similar to (\ref{eq:CovarianceFOEMWhiteNoise}) in combination with Fourier transforms of radial functions. The spatial behavior of this model can be described following equation (\ref{eq:SPDESpatialTraceGenOEM}), with spatial White Noise source term, $X_{S} = W_{S}$. In general, this Langevin equation model is not an evolving Matérn model. It is the particular case when the parameter $\nu$, called \textit{curvature coefficient}, equals to $0$.

\begin{Example}{\textbf{Heat equation.}}
\label{Ex:HeatEquation}
\end{Example}

We now consider the stochastic Heat (or Diffusion) Equation over $\mathbb{R}^{d}\times\mathbb{R}$
\begin{equation}
\label{eq:SPDEHeatEquation}
\dfrac{\partial U }{\partial t} - a \Delta U = X,
\end{equation}
where $a > 0 $ is the \textit{diffusivity} parameter. It is a first order evolution model with spatial symbol $g(\xi) = a|\xi|^{2}$. In this case, the spatio-temporal symbol function $(\xi , \omega) \mapsto i \omega + a|\xi|^{2}$ is not strictly positive, the origin being the only zero of $g$. There is thus no uniqueness of stationary solutions, if they exist. Using similar arguments as those used in Section \ref{Section:MaternWithoutKappa}, one can see that the only stationary solutions to the homogeneous Heat Equation 
\begin{equation}
\label{eq:SPDEHeatEquationHomogeneous}
\dfrac{\partial U_{H} }{\partial t} - a \Delta U_{H} = 0
\end{equation}
are random constants. Because of the singularity at the origin of the function $| g |^{-2}$, the existence condition (\ref{eq:ConditionExistence}) does not always hold. Existence needs to be checked for each source term $X$. Let us first consider that the source term is a spatio-temporal White Noise. Equation (\ref{eq:SPDEHeatEquation}) becomes 
\begin{equation}
\label{eq:SPDEHeatEquationWhiteNoise}
\dfrac{\partial U }{\partial t} - a \Delta U = W.
\end{equation}
Using Theorem \ref{TheoremExistAndUniqueStationary}, one concludes (see Appendix \ref{Section:ExistenceSolutionsHeatEquationWhiteNoise}) that \textit{there exists stationary solutions to the stochastic Heat equation (\ref{eq:SPDEHeatEquationWhiteNoise}) only for spatial dimensions $d \geq 3 $}, and in those cases, they can only be conceived  as GeRFs and never as random functions continuous in mean-square.
When $d=3$, computations reported in Appendix \ref{Section:ObtainingCovarianceHeatEquationWhiteNoised3} show that the covariance structure  is
\begin{equation}
\label{eq:CovarianceSPDEHeatEquationWhiteNoise}
\rho_{U}^{W}(h,u)= \dfrac{1}{(2\pi)^{\frac{d+1}{2}}} \dfrac{\pi}{2a|h|}\erf\left( \frac{|h|}{2\sqrt{a|u|}} \right).
\end{equation}
This covariance must be interpreted in a distributional sense, since it is not defined at $|h| = |u| = 0 $. A spatial trace of the stationary field associated to (\ref{eq:CovarianceSPDEHeatEquationWhiteNoise}), $U_{S}$, can be described evaluating this covariance in $u=0$ with $h \neq 0$. We obtain that $U_{S}$ satisfies the spatial SPDE
\begin{equation}
\sqrt{2}(-\Delta)^{\frac{1}{2}}U_{S} \stackrel{2nd \ o.}{=} W_{S},
\end{equation}
where $ W_{S}$ is a spatial White Noise. In other words, $U_{S}$ is a Matérn model without range parameter as presented in Section \ref{Section:MaternWithoutKappa}. See Eq. (\ref{eq:CovarianceMaternWithoutKappa}).

When $X$ is an arbitrary source term, Theorem \ref{Theorem:CovarianceCovolutionWhiteNoiseSourceTerm} cannot be applied for spatial dimensions smaller that $3$. For $d=3$, a convolvability condition between $\rho_{X}$ and (\ref{eq:CovarianceSPDEHeatEquationWhiteNoise}) must be satisfied. Nevertheless, the existence of a solution can be ensured independently on existence of solutions with White Noise source term by imposing necessary conditions on $\mu_{X}$ such that the existence criteria (\ref{eq:ConditionExistence}) holds. For example, one could require $\mu_{X}$ to be null in some neighborhood of the origin.

\begin{Example}{\textbf{Wave equation.}}
\label{Ex:WaveEquation}
\end{Example}

As a final example we consider the stochastic wave equation
\begin{equation}
\label{eq:SPDEWaveEquation}
\dfrac{\partial^{2} U}{\partial t^{2}} - c^{2}\Delta U  = X,
\end{equation}
where $X$ is a stationary random field and $c > 0$ is the \textit{propagation velocity}. This is a second order evolution model with spatial symbol function $g(\xi) = c^{2}|\xi|^{2} $. The null-set of the associated spatio-temporal symbol function $(\xi , \omega) \mapsto -\omega^{2} + c^{2}|\xi|^{2}$ is the spatio-temporal cone $\mathcal{C} = \lbrace (\xi,\omega) \in \mathbb{R}^{d}\times\mathbb{R} \mid  |\omega| = c |\xi|  \rbrace  $. As a consequence, uniqueness of a potential stationary solution does not hold. Following Remark \ref{Rem:gWithZerosHomogeneousProblem}, stationary solutions to the homogeneous wave equation
\begin{equation}
\label{eq:SPDEWaveEquationHomogeneous}
\dfrac{\partial^{2} U_{H}}{\partial t^{2}} - c^{2}\Delta U_{H}  = 0
\end{equation}
are found by studying covariance structures associated to spectral measures supported on the cone $\mathcal{C}$. A spectral measure $\mu_{U_{H}}$ over $\mathbb{R}^{d}\times\mathbb{R}$ supported on $\mathcal{C}$, can be described trough its action over test functions $\psi \in \mathcal{S}(\mathbb{R}^{d}\times\mathbb{R})$ by
\begin{equation}
\label{eq:MeasureSupportedCone}
\langle \mu_{U_{H}} , \psi \rangle = \sqrt{2\pi}\int_{\mathbb{R}^{d}}\frac{\psi(\xi , c|\xi|) + \psi(\xi , -c| \xi |)}{2} d \mu_{U_{H}^{S}}(\xi) = \sqrt{2\pi}\int_{\mathbb{R}^{d}}\int_{\mathbb{R}} \psi(\xi , \omega) d \left( \frac{\delta_{-c|\xi|} + \delta_{c|\xi|}}{2} \right)(\omega)d \mu_{U_{H}^{S}}(\xi),
\end{equation}
where $\mu_{U_{H}^{S}}$ is a spectral measure over $\mathbb{R}^{d}$. In the right hand side of (\ref{eq:MeasureSupportedCone}), the disintegration language is used for the measure $\mu_{U_{H}}$,  which can then be expressed as 
\begin{equation}
d\mu_{U_{H}}(\xi , \omega) = \sqrt{2\pi}d\left(\frac{\delta_{-c|\xi|} + \delta_{c|\xi|}}{2}\right)(\omega)d \mu_{U_{H}^{S}}(\xi).
\end{equation}
Hence, all stationary solutions of (\ref{eq:SPDEWaveEquationHomogeneous}) have a spectral measure of this form.  After applying a temporal Fourier transform, the associated covariance is 
\begin{equation}
\label{eq:CovarianceHomogeneousWaveEquation}
\rho_{U_{H}}(h,u) = \mathscr{F}_{S}\left( \xi \mapsto \cos(c|\xi||u|)d\mu_{U_{H}^{S}}(\xi)  \right)(h).
\end{equation}
The factor $\sqrt{2\pi}$ is introduced in (\ref{eq:MeasureSupportedCone}) in order to identify $\mu_{U_{H}^{S}}$  as the spectral measure of a spatial trace of $U_{H}$, as it can be verified by evaluating (\ref{eq:CovarianceHomogeneousWaveEquation}) at $u=0$. $\mu_{U_{H}^{S}}$ can be chosen arbitrarily. Thus, we can use any spatial stationary model to construct a spatio-temporal stationary solution of (\ref{eq:SPDEWaveEquationHomogeneous}) following this spatial model for any fixed time coordinate.
As an example, we call \textit{waving Matérn model} the model with spatial spectral measure,  $d\mu_{U_{H}^{S}}(\xi) = (2\pi)^{-d/2}a(\kappa^{2} + |\xi|^{2})^{-\alpha}d\xi  $, with $a , \kappa > 0 $ and $\alpha \in \mathbb{R} $. The associated covariance is
\begin{equation}
\label{eq:CovarianceWavingMaternModel}
\rho(h,u) = \mathscr{F}_{S}\left( \xi \mapsto \dfrac{\cos( c|\xi||u| )}{ (2\pi)^{\frac{d}{2}}a(\kappa^{2} + |\xi|^{2})^{\alpha} }    \right)(h).
\end{equation}
Let us now go back to the existence of stationary solutions of (\ref{eq:SPDEWaveEquation}) with $X=W$, i.e.
\begin{equation}
\label{eq:SPDEWaveEquationWhiteNoise}
\dfrac{\partial^{2} U}{\partial t^{2}} - c^{2}\Delta U  = W.
\end{equation}
Since the function $(\xi,\omega) \mapsto (-\omega^{2} + c^{2}|\xi|^{2})^{-2}$ is not locally integrable, by applying Theorem \ref{TheoremExistAndUniqueStationary} we conclude that \textit{there are no stationary solutions to the stochastic wave equation (\ref{eq:SPDEWaveEquationWhiteNoise})}. Hence, we cannot apply Theorem \ref{Theorem:CovarianceCovolutionWhiteNoiseSourceTerm} to relate the covariance of a possible stationary solution of (\ref{eq:SPDEWaveEquation}) to the covariance of the solution with White Noise source term. The existence of a stationary solution to (\ref{eq:SPDEWaveEquation}) must be then studied for every particular case of $X$. Notice however that the existence is guaranteed when the supports of the spectral measure of the source term and the spatio-temporal cone $\mathcal{C}$ are separated by neighborhoods.

\section{Conclusion}
\label{Section:Conclusion}

We have proposed a very general setting that allows to relate a SPDE to spatial and spatio-temporal covariance structures through the specification of symbol functions. It is grounded on the concept of Generalized Random Field, stochastic analogue to Schwartz's concept of distribution \cite{Schwartz1966} as already proposed in \citet{Ito1954} and \citet{Matheron1965}. This setting offers a convenient framework to build and characterize models of random fields that are stationary solutions, when they exist, of a very large class of SPDEs. Their covariance structure is in direct relationship with the symbol function thanks to Theorem \ref{TheoremExistAndUniqueStationary}. In particular, this setting allows to handle relatively easily SPDEs with fractional behavior, in time, in space, and in both spatial and temporal dimensions. Thanks to this framework, we were able to construct very general models, that include and encompass existing models, as shown in details in Section \ref{Section:IllustrationsSpatial} and Section \ref{Section:SpatioTemporalModels}. 

Theorem \ref{Theorem:CovarianceCovolutionWhiteNoiseSourceTerm} establishes that, under mild conditions, the covariance of the stationary solution of a given SPDE for general random source term with covariance $\rho_X$ is the convolution between the covariance of the same SPDE with White Noise source term and $\rho_X$. This results is a powerful tool for easily characterizing solutions of very general SPDEs. It also emphasizes the central role played by a White Noise source term.

We envision this work as a contribution strengthening the SPDE paradigm shift for analyzing spatial and spatio-temporal data as initiated in \cite{Lindgren2011}. Our contribution offers the possibility to build and characterize models far beyond the Matérn family which is currently the  covariance model considered within most SPDE implementations.

Efficient simulation of our models can be easily conceived using Fourier analysis based PDE-solvers as proposed in \cite{Lang2011}. Inference and simulation methods presented in \citet{Sigrist2015} can be easily adapted to any first order evolution models presented in Section \ref{Section:EvolutionModels}. Since the linear operators considered in this work are not strictly speaking differential operators, methods inspired by the Finite Elements Method or by the Finite Difference Method are not applicable without specific adaptation. For instance, \citet{Bolin2017} propose adaptations of the Finite Elements Method for Matérn models with fractional regularity. 

\begin{appendices}

\section{Reminders on tempered distributions}
\label{Section:TemperedDistribution}

Here we give a brief overview of the main definitions and results regarding Schwartz's theory of Distributions in a tempered framework. For a more detailed presentation, the reader is referred to \cite{Donoghue1969} and, of course, to \cite{Schwartz1966}. For a brief introduction with geostatistical purposes, we suggest \citet[Appendix A]{Matheron1965}.

Let $\mathcal{S}(\mathbb{R}^{d})$ be the set of all complex, smooth and fast decreasing  functions over $\mathbb{R}^{d}$,
$$\mathcal{S}(\mathbb{R}^{d}) = \lbrace \varphi \in C^{\infty}(\mathbb{R}^{d}) \quad \big| \quad \| x^{\alpha}D^{\beta}\varphi \|_{\infty} = \displaystyle\sup_{x \in \mathbb{R}^{d}}|x^{\alpha}D^{\beta}\varphi(x)| < \infty,\ \forall \alpha , \beta \in \mathbb{N}^{d} \rbrace,$$ 
where the multi-index notation for the power $x^{\alpha}$ and the differential operator $D^{\beta}$ for $\alpha, \beta \in \mathbb{N}^{d}$ is used, meaning respectively $ x^{\alpha} = x_{1}^{\alpha_{1}} \cdots x_{d}^{\alpha_{d}}$ and $D^{\beta} = \frac{ \partial^{|\beta|}}{ \partial x^{\beta_{1}}_{1} \cdots \partial x^{\beta_{d}}_{d}}$, with $|\beta| = \beta_{1} + ... + \beta_{d}$. Equipped with a suitable topology, $\mathcal{S}(\mathbb{R}^{d})$ is a complete metric space, known as the Schwartz space of \textit{test functions}. Its dual space, i.e. the space of all continuous linear functionals from $\mathcal{S}(\mathbb{R}^{d})$ to $\mathbb{C}$, is called the space of \textit{tempered distributions} and it is denoted by $\mathcal{S}'(\mathbb{R}^{d})$. In order to emphasize the dual aspect of tempered distributions and test functions, we will denote $\langle T,\varphi \rangle$ the action of $T \in {\cal S}'(\mathbb{R}^{d})$ on $\varphi \in {\cal S}(\mathbb{R}^{d})$.

Tempered distributions can be seen as a generalization of functions, on which the Fourier transform and differentiations of any order can be rigorously defined. Polynomials, continuous and bounded functions or functions $f \in L^{p}(\mathbb{R}^{d})$ with $p \in [ 1 , \infty ]$ can be interpreted as tempered distributions through the integral $\langle f , \varphi \rangle := \int_{\mathbb{R}^{d}}f(x)\varphi(x)dx $ which is well-defined for all $\varphi \in \mathcal{S}(\mathbb{R}^{d})$. Similarly, a finite measure $\mu$ over $\mathbb{R}^{d}$ can also be interpreted as a tempered distribution through the integral $\langle \mu , \varphi \rangle := \int_{\mathbb{R}^{d}}\varphi(x)d\mu(x) $,  $\varphi \in \mathcal{S}(\mathbb{R}^{d})$.

Tempered distributions can be differentiated any number of times. Let $D^{\alpha}$ be a differential operator with $\alpha \in \mathbb{N}^{d}$. Inspired by the integration by parts formula, the derivative of a tempered distribution $T \in \mathcal{S}'(\mathbb{R}^{d})$ is defined as a new tempered distribution $D^{\alpha}T \in \mathcal{S}'(\mathbb{R}^{d})$  through $\langle D^{\alpha}T , \varphi \rangle := (-1)^{|\alpha|}\langle T , D^{\alpha}\varphi \rangle$ for all $\varphi \in \mathcal{S}(\mathbb{R}^{d})$. The Fourier transform and its inverse are defined for any test function $\varphi \in \mathcal{S}(\mathbb{R}^{d})$ as
\begin{equation}
\label{eq:DefFourierTestFunction}
\mathscr{F}(\varphi)(\xi) = \dfrac{1}{(2\pi)^{d/2}}\int_{\mathbb{R}^{d}}e^{-i\xi^{T}x}\varphi(x)dx,  \qquad \mathscr{F}^{-1}(\varphi)(\xi) = \dfrac{1}{(2\pi)^{d/2}}\int_{\mathbb{R}^{d}}e^{i\xi^{T}x}\varphi(x)dx.
\end{equation}
For tempered distributions, the Fourier transform is defined as a new tempered distribution through the transfer formula
\begin{equation}
\langle \mathscr{F}(T) , \varphi \rangle :=  \langle T , \mathscr{F}(\varphi) \rangle; \quad  \langle \mathscr{F}^{-1}(T) , \varphi \rangle :=  \langle T , \mathscr{F}^{-1}(\varphi) \rangle, 
\quad \forall \varphi \in \mathcal{S}(\mathbb{R}^{d}),\ T \in \mathcal{S}'(\mathbb{R}^{d}).
\end{equation}
The Fourier transform is a continuous bijective endomorphism over $\mathcal{S}(\mathbb{R}^{d})$ and over $\mathcal{S}'(\mathbb{R}^{d})$. The classical property of the Fourier transform with respect to the differentiation, $\mathscr{F}( D^{\alpha} T ) = i^{|\alpha |}\xi^{\alpha}\mathscr{F}(T)$, where $\xi$ denotes the variable in the frequency space, holds also for every tempered distribution $T$.

Let us also define the space $\mathcal{O}_{M}(\mathbb{R}^{d})$ of complex smooth functions defined over $\mathbb{R}^{d}$ such that all of its derivatives are polynomially bounded. Explicitly,
$$\mathcal{O}_{M}(\mathbb{R}^{d}) = \lbrace f \in C^{\infty}(\mathbb{R}^{d}) \ : \ \forall \alpha \in \mathbb{N}^{d} \ \exists C > 0 \ \exists N \in \mathbb{N} \ \hbox{such that} \  |D^{\alpha}f(x)|\leq C(1+|x|^{2})^{N} \ \forall x \in \mathbb{R}^{d}  \rbrace.  $$
This space is known as  \textit{the space of multiplicators} of the Schwartz space. If $f \in \mathcal{O}_{M}(\mathbb{R}^{d})$ and $\varphi \in \mathcal{S}(\mathbb{R}^{d})$, then $f\varphi \in \mathcal{S}(\mathbb{R}^{d})$. If $T \in \mathcal{S}'(\mathbb{R}^{d})$, the product $fT \in \mathcal{S}'(\mathbb{R}^{d})$ is defined through $\langle fT , \varphi \rangle = \langle T , f\varphi \rangle $ for every $\varphi \in \mathcal{S}(\mathbb{R}^{d})$. If $ f \in \mathcal{O}_{M}(\mathbb{R}^{d})$, then its Fourier transform $\mathscr{F}(f)$ is convolvable with any tempered distribution, and the exchange formula for the Fourier transform holds: $\mathscr{F}(fT) =  (2\pi)^{-\frac{d}{2}}\mathscr{F}(f)\ast\mathscr{F}(T)$ for every $T \in \mathcal{S}'(\mathbb{R}^{d})$. See \cite{Schwartz1966}, chapter VII, section $5$ and Theorem $XV$ in section $8$. 

\section{Proof of Proposition \ref{Prop:MultiplicationgZ} }
\label{Section:ProofPropositiongZ}

The main difficulty of this Proposition lies in a proper definition of the product $gZ$ as a GeRF. Indeed, we could simply write $\langle gZ , \varphi \rangle := \langle Z , g\varphi \rangle $, but $Z$ is only defined over functions in $\mathcal{S}(\mathbb{R}^{d})$, and $g\varphi$ is not in general in $ \mathcal{S}(\mathbb{R}^{d})$. Nevertheless, we will show that we can define $ \langle Z , f \rangle  $ if $Z$ is a slow-growing random measure and $f$ is a continuous function with fast decreasing behavior .

We define $C_{FD}(\mathbb{R}^{d}) := \lbrace f \in C(\mathbb{R}^{d}) \ \ | \ \| (1+|x|^{2})^{N}f \|_{\infty} < \infty \ \forall N \in \mathbb{N} \rbrace $, the space of all continuous functions with fast decreasing behavior, equipped with the following topology: a sequence of functions $(f_{n})_{n \in \mathbb{N}} \subset C_{FD}(\mathbb{R}^{d}) $ converges to $f \in C_{FD}(\mathbb{R}^{d})$, denoted $ f_{n} \stackrel{C_{FD}}{\to} f $, if for all $N \in \mathbb{N}$ we have that $ \| (1 + |x|^{2})^{N}(f_{n}-f) \|_{\infty} \to 0 $. This topology is induced by the metric
\begin{equation}
(f,g) \mapsto \sum_{N \in \mathbb{N}} \frac{1}{2^{N}} \frac{\| (1 + |x|^{2})^{N}(f-g) \|_{\infty}}{1 + \| (1 + |x|^{2})^{N}(f-g) \|_{\infty} }, \ \quad f,g \in C_{FD}(\mathbb{R}^{d}).
\end{equation}
For this topological vector space, the following two lemmas hold (proofs given below). 
\begin{Lemma}
\label{LemmaSdenseCdr}
$\mathcal{S}(\mathbb{R}^{d}) \subset C_{FD}(\mathbb{R}^{d}) $, and it is a dense sub-space (with the topology of $C_{FD}$).
\end{Lemma}
\begin{Lemma}
\label{LemmaMuContinuousFunctional}
$\mathcal{M}_{SG}(\mathbb{R}^{d}) = C_{FD}'(\mathbb{R}^{d})$, that is, every measure $\mu \in \mathcal{M}_{SG}(\mathbb{R}^{d})$ defines a continuous linear functional $T$ over  $C_{FD}(\mathbb{R}^{d})$ through the integral 
\begin{equation}
\label{eq:SlowGrowingMeasureDualFastDecreasing}
\langle T , f \rangle = \int_{\mathbb{R}^{d}}f(x)d\mu(x), \quad \forall f \in C_{FD}(\mathbb{R}^{d}).
\end{equation}
Conversely, for every continuous linear functional $T : C_{FD}(\mathbb{R}^{d}) \to \mathbb{C}$ there exists a unique $\mu \in \mathcal{M}_{SG}(\mathbb{R}^{d})$ such that (\ref{eq:SlowGrowingMeasureDualFastDecreasing}) holds.
\end{Lemma}

We now prove Proposition  \ref{Prop:MultiplicationgZ}. If $g$ is a continuous function bounded by a polynomial, then $g\varphi \in C_{FD}(\mathbb{R}^{d})$ for all $\varphi \in \mathcal{S}(\mathbb{R}^{d})$. Since, as stated in Lemma \ref{LemmaSdenseCdr}, $\mathcal{S}(\mathbb{R}^{d})$ is dense in $C_{FD}(\mathbb{R}^{d})$ , we can construct the random variable $\langle Z , g\varphi \rangle$ as a limit in a mean-square sense. Let $(g_{n})_{n \in \mathbb{N}} \subset \mathcal{S}(\mathbb{R}^{d})$ be a sequence such that $g_{n} \stackrel{C_{FD}}{\to} g\varphi $. Consider the sequence of square-integrable (centered) random variables $(\langle Z , g_{n} \rangle)_{n \in \mathbb{N}}$. We obtain by linearity that 
\begin{equation}
\label{eq:Z(gn)Cauchy}
\mathbb{E}( \left| \langle Z , g_{n} \rangle - \langle Z , g_{m} \rangle \right|^{2}  ) = \mathbb{E}( \left| \langle Z , g_{n} - g_{m} \rangle \right|^{2} ) = \int_{\mathbb{R}^{d}\times\mathbb{R}^{d}}(g_{n}-g_{m})(x)\overline{(g_{n}-g_{m})}(y)dC_{Z}(x,y).
\end{equation} 
Since the sequence  $(g_{n})_{n \in \mathbb{N}}$ converges in $C_{FD}(\mathbb{R}^{d})$, it is tedious but easy to show that the sequence $(g_{n}\otimes\overline{g_{n}})_{n \in \mathbb{N} }$ converges in $C_{FD}(\mathbb{R}^{d}\times\mathbb{R}^{d})$. Since $C_{Z} \in \mathcal{M}_{SG}(\mathbb{R}^{d}\times\mathbb{R}^{d})$, by Lemma \ref{LemmaMuContinuousFunctional} the integral in (\ref{eq:Z(gn)Cauchy}) goes to zero as $n$ and $m$ grow, due to the continuity of $C_{Z}$ interpreted as a linear functional over $C_{FD}(\mathbb{R}^{d}\times\mathbb{R}^{d})$. The sequence $(\langle Z , g_{n} \rangle)_{n \in \mathbb{N}}$ is thus a Cauchy sequence in $L^{2}(\Omega , \mathcal{F} , \mathbb{P})$. Hence, it is convergent, and we write 
\begin{equation}
\label{eq:DefgZ}
\langle gZ , \varphi \rangle := \langle Z , g\varphi \rangle := \lim_{n \to \infty} \langle Z , g_{n} \rangle,
\end{equation}
where the limit is taken in the $L^{2}$ sense. This limit does not depend on the sequence converging to $g\varphi$ chosen. Indeed, if $(f_{n})_{n \in \mathbb{N}} \subset \mathcal{S}(\mathbb{R}^{d})$ is another sequence such that $f_{n} \stackrel{C_{FD}}{\to} g\varphi $, then the sequence defined by $h_{n} = f_{n} - g_{n}$ for all $n \in \mathbb{N}$ converges to $0$ in $C_{FD}(\mathbb{R}^{d})$. As well, the sequence $ h_{n} \otimes \overline{h_{n}}$ converges to $0$ in $C_{FD}(\mathbb{R}^{d} \times \mathbb{R}^{d})$, and by the same argument used to prove the convergence to $0$ of (\ref{eq:Z(gn)Cauchy}), we have $\langle Z , h_{n} \rangle  \to 0$ in the $L^{2}$ sense. Thus $\lim_{n \to \infty} \langle Z , f_{n} \rangle = \lim_{n \to \infty} \langle Z , g_{n} \rangle = \langle Z , g\varphi \rangle$. 

The covariance structure of $gZ$ can be easily obtained as a limit of covariances. Thus we obtain,
\begin{equation}
\label{CovgZ}
\mathbb{C}ov( \langle gZ , \varphi \rangle , \langle gZ , \phi \rangle  ) = \int_{\mathbb{R}^{d}\times\mathbb{R}^{d}}\varphi(x)\overline{\phi}(y)g(x)\overline{g}(y)dC_{Z}(x,y) = \langle (g\otimes\overline{g})C_{Z} , \varphi \otimes \overline{\phi} \rangle, \quad \forall \varphi, \phi \in \mathcal{S}(\mathbb{R}^{d}).
\end{equation}
The result of Corollary \ref{CorolLemmaMultiplicationgZ}, which describes the case of a slow-growing orthogonal random measure, follows directly from (\ref{CovgZ}). Details are left to the reader. $\blacksquare$ 

\medskip
\textbf{Proof of Lemma \ref{LemmaSdenseCdr}.}
It is clear that $\mathcal{S}(\mathbb{R}^{d}) \subset C_{FD}(\mathbb{R}^{d})$. To prove the density, we first prove that if $f \in C_{FD}(\mathbb{R}^{d})$ and $\varphi \in \mathcal{S}(\mathbb{R}^{d})$, then $ f \ast \varphi \in \mathcal{S}(\mathbb{R}^{d})$. It is clear that $f$ is integrable and bounded, as well as $\varphi$ which, in addition, is smooth. Thus $f \ast \varphi $ is a smooth integrable and bounded function, and its Fourier transform satisfies $\mathscr{F}( f \ast \varphi) = (2\pi)^{\frac{d}{2}}\mathscr{F}(f)\mathscr{F}(\varphi)$. We have that $\mathscr{F}(\varphi) \in \mathcal{S}(\mathbb{R}^{d})$ since $\mathscr{F}$ is a bijective endomorphism  of $\mathcal{S}(\mathbb{R}^{d})$. Since $f \in C_{FD}(\mathbb{R}^{d})$ we conclude by Riemann-Lebesgue lemma that $\mathscr{F}( f )$ is a smooth function with all its derivatives vanishing at infinity. Thus $\mathscr{F}(f) \in \mathcal{O}_{M}(\mathbb{R}^{d})  $, which implies that $(2\pi)^{\frac{d}{2}}\mathscr{F}(f)\mathscr{F}(\varphi) \in \mathcal{S}(\mathbb{R}^{d})$. This proves that $f \ast \varphi = \mathscr{F}^{-1}\left(  (2\pi)^{\frac{d}{2}}\mathscr{F}(f)\mathscr{F}(\varphi) \right) \in \mathcal{S}(\mathbb{R}^{d})$.

Let $(\phi_{n})_{n \in \mathbb{N}} \subset \mathcal{S}(\mathbb{R}^{d})$ be a regularizing sequence of positive compactly supported smooth functions, such that $\supp(\phi_{n}) = \overline{B_{1/n}(0)}$ and  $\int_{\mathbb{R}^{d}}\phi_{n}(x)dx = 1 $ for all $n\in\mathbb{N}$. Here $B_{r}(0) \subset \mathbb{R}^{d}$ denotes the open ball with center $0$ and radius $r > 0$. We consider the sequence of functions $f_{n} = f \ast \phi_{n}$, which are all in $\mathcal{S}(\mathbb{R}^{d})$. We will prove that $f_{n} \stackrel{C_{FD}}{\to} f $. Let $m \in \mathbb{N}$ be fixed. We must show that $\| (1 + |x|^{2})^{m}(f_{n} - f) \|_{\infty} \to 0 $ as $n \to \infty $. Let $\epsilon > 0 $. As $f \in  C_{FD}(\mathbb{R}^{d})$, we can take $R > 0 $ large enough such that for every $x$ such that $|x| > R-1$, $(1+2|x|^{2})^{m}|f(x)| < \frac{\epsilon}{3(2^{m-1}+2^{2m-1})} $ holds. Notice that in this case, $(1+|x|^{2})^{m}|f(x)| < \epsilon/3 $. Since $f$ is continuous, it is uniformly continuous over the compact set $\overline{B_{R+1}(0)}$.  Thus, there exists $\delta > 0 $ such that if $|x-y|<\delta$, then $|f(x)-f(y)|< \frac{\epsilon}{3(1+R^{2})^{m}}$ for all $x,y \in \overline{B_{R+1}(0)}$. Consider $n_{0} \in \mathbb{N}$ such that $1/n_0 < \delta $. Then,  for all $n \geq n_{0}$,
\begin{align}
\| (1+|x|^{2})^{m}(f - f_{n}) \|_{\infty} & = & \sup_{x \in \mathbb{R}^{d}}\left| \int_{B_{1/n}(0)}(1+|x|^{2})^{m}(f(x)-f(x-y))\phi_{n}(y)dy \right| \nonumber  \\ 
 & \leq &  \sup_{x \in \overline{B_{R}(0)}}\left| \int_{B_{1/n}(0)}(1+|x|^{2})^{m}(f(x)-f(x-y))\phi_{n}(y)dy \right| \nonumber \label{PartA} \tag{a}\\
  &  & + \sup_{x \in \overline{B_{R}(0)}^{c}}\left| \int_{B_{1/n}(0)}(1+|x|^{2})^{m}(f(x)-f(x-y))\phi_{n}(y)dy \right|. \nonumber \label{PartB} \tag{b}  \\ \label{EqsBoundingUnifConv}
\end{align}
For the first term (\ref{PartA})  uniform continuity of $f$ implies
\begin{equation}
\label{EqBoundingPartA}
\sup_{x \in \overline{B_{R}(0)}}\left| \int_{B_{1/n}(0)}(1+|x|^{2})^{m}(f(x)-f(x-y))\phi_{n}(y)dy \right| \leq \int_{B_{1/n}(0)}(1+R^{2})^{m}\dfrac{\epsilon}{3(1+R^{2})^{m}}\phi_{n}(y)dy = \dfrac{\epsilon}{3}. 
\end{equation}
Regarding the second term (\ref{PartB}), the integral is split  to obtain
\begin{equation}
\label{EqBoundingPartB}
(\ref{PartB}) \leq \sup_{x \in \overline{B_{R}(0)}^{c}}\Big\lbrace \underbrace{\int_{B_{1/n}(0)}(1+|x|^{2})^{m}|f(x)|\phi_{n}(y)dy}_{\leq \frac{\epsilon}{3} } + \int_{B_{1/n}(0)}(1+|x|^{2})^{m}|f(x-y)|\phi_{n}(y)dy  \Big\rbrace
\end{equation}
When applying Jensen's inequality twice, one shows that $(1+|x|^{2})^{m} \leq 2^{m-1}[ (1+2|x-y|^{2})^{m} + 2^{m}|y|^{2m} ] $ for all $x$ and $y$, and thus
\begin{align}
\int_{B_{1/n}(0)}(1+|x|^{2})^{m}|f(x-y)|\phi_{n}(y)dy & \leq  & 2^{m-1}\biggl[ \int_{B_{1/n}(0)}\underbrace{(1+2|x-y|^{2})^{m}|f(x-y)|}_{< \frac{\epsilon}{3(2^{m-1}+2^{2m-1})} \text{ from } |x-y| > R-1 }\phi_{n}(y)dy  \nonumber \\
 &   &  +  2^{m}\int_{B_{1/n}(0)}\underbrace{|f(x-y)|}_{ < \frac{\epsilon}{3(2^{m-1}+2^{2m-1})}}\underbrace{|y|^{2m}}_{\leq 1}\phi_{n}(y)dy \biggr] \nonumber \\
   &  <  & 2^{m-1}\left(\dfrac{\epsilon}{3(2^{m-1}+2^{2m-1})} + 2^{m}\dfrac{\epsilon}{3(2^{m-1}+2^{2m-1})}\right) = \dfrac{\epsilon}{3}. \nonumber \\ \label{EqsBoundingIntegrals}
\end{align}
Hence considering (\ref{EqBoundingPartB}) and (\ref{EqsBoundingIntegrals}) we finally obtain $(\ref{PartB}) < 2\epsilon/3 $. Putting together this result and (\ref{EqBoundingPartA}) on equation (\ref{EqsBoundingUnifConv}), we finally obtain  that for all $n \geq n_{0}$,
\begin{equation}
\| (1+|x|^{2})^{m}(f - f_{n}) \|_{\infty} < \epsilon,
\end{equation}
hence $ \| (1+|x|^{2})^{m}(f - f_{n}) \|_{\infty} \to 0 $. Since $m$ was arbitrary, this result holds for all $m$. We therefore conclude that $f_{n} \stackrel{C_{FD}}{\to} f $. Since for any arbitrary $f \in C_{FD}(\mathbb{R}^{d})$ we can find a sequence included in $\mathcal{S}(\mathbb{R}^{d})$ which converges to $f$, we conclude that $\mathcal{S}(\mathbb{R}^{d})$ is dense in $C_{FD}(\mathbb{R}^{d})$.
$\blacksquare$

\textbf{Proof of Lemma \ref{LemmaMuContinuousFunctional}.} This Lemma is an analogue of the famous Riesz's Representation Theorem, where the duals of some spaces of continuous function are identified with some spaces of measures. This result is probably trivial for an analyst, but we have not found any reference where the way to obtain it is explicited.

We first make the following claim: if $T : C_{FD}(\mathbb{R}^{d}) \to \mathbb{C}$ is a linear functional, then it is continuous if and only if there exists $C>0$ and $N_0 \in \mathbb{N}$ such that
\begin{equation}
\label{EqPropDescriptionContinuousFunctionals}
| \langle T , f \rangle | \leq C \| (1+|x|^{2})^{N_0}f \|_{\infty} \quad \forall f \in C_{FD}(\mathbb{R}^{d}).
\end{equation}
The sufficiency of this claim is straightforward. Indeed, let consider any sequence $(f_{n})_{n \in \mathbb{N}}$ such that $ f_{n} \stackrel{C_{FD}}{\to} 0 $. In particular this sequence satisfies $\| (1+|x|^{2})^{N_0}f_{n} \|_{\infty} \to 0 $, thus $ \langle T , f_{n} \rangle \to 0 $, and $T$ is continuous. Let us prove the necessity. We argue by contradiction. Let us suppose that $T$ is continuous but that for all $C >0$ and for all $N \in \mathbb{N}$ we can find a function $ f_{C,N} \in C_{FD}(\mathbb{R}^{d}) $ such that $|\langle T , f_{C,N} \rangle | > C \| (1+|x|^{2})^{N}f_{C,N} \|_{\infty} $. We consider $C = n^{2}$ and $N=n$ for all $n \in \mathbb{N}$. We obtain thus a sequence of functions $(f_{n})_{n \in \mathbb{N}}$ in $C_{FD}(\mathbb{R}^{d})$ such that
\begin{equation}
\label{EqTf>n^2|(1+x2)^nfn|}
|\langle T , f_{n} \rangle | > n^{2} \| (1+|x|^{2})^{n}f_{n} \|_{\infty} \quad \forall n \in \mathbb{N}.
\end{equation}
Let us define the sequence of functions $(\phi_{n})_{n \in \mathbb{N}}$ by
\begin{equation}
\label{EqDefPhinProofProp}
\phi_{n} = \frac{f_{n}}{n\displaystyle\sum_{m \leq n}\| (1+|x|^{2})^{m}f_{n} \|_{\infty}}, \quad n \in \mathbb{N}.
\end{equation}
Clearly the sequence $(\phi_{n})_{n \in \mathbb{N}}$ is in $C_{FD}(\mathbb{R}^{d})$. Let $M \in \mathbb{N}$. By (\ref{EqDefPhinProofProp}), if $n \geq M$ it holds that
\begin{equation}
\| (1+|x|^{2})^{M}\phi_{n} \|_{\infty} = \frac{ \| (1+|x|^{2})^{M}f_{n} \|_{\infty} }{n\displaystyle\sum_{m \leq n}\| (1+|x|^{2})^{m}f_{n} \|_{\infty}} < \frac{1}{n},
\end{equation}
and thus $ \phi_{n} \stackrel{C_{FD}}{\to} 0 $. On the other hand, by (\ref{EqTf>n^2|(1+x2)^nfn|}) we get 
\begin{equation}
| \langle T , \phi_{n} \rangle | = \dfrac{1}{n}\dfrac{| \langle T , f_{n} \rangle |}{\displaystyle\sum_{m \leq n}\| (1+|x|^{2})^{m}f_{n} \|_{\infty}} > \dfrac{n^{2}}{n}\dfrac{\| (1+|x|^{2})^{n}f_{n} \|_{\infty}}{\displaystyle\sum_{m \leq n}\| (1+|x|^{2})^{m}f_{n} \|_{\infty}} \geq 1,
\end{equation}
where we have used that 
\begin{equation}
\displaystyle\sum_{m \leq n} \| (1+|x|^{2})^{m}f_{n} \|_{\infty} \leq n\| (1+|x|^{2})^{n}f_{n} \|_{\infty} 
\end{equation}
since $ \| (1+|x|^{2})^{m}f_{n} \|_{\infty} \leq \| (1+|x|^{2})^{n}f_{n} \|_{\infty} $ for all $m \leq n$. We conclude that $| \langle T , \phi_{n} \rangle |$ does not converge to $0$ as $n$ grows, and thus $T$ is not continuous, which is a contradiction. Hence, our claim holds.

Let us now prove the Lemma. Let $\mu \in \mathcal{M}_{SG}(\mathbb{R}^{d})$. By definition, there exists $N \in \mathbb{N}$ such that $(1+|x|^{2})^{-N}|\mu| $ is finite. By taking $C = \left( (1+|x|^{2})^{-N}|\mu| \right)(\mathbb{R}^{d}) < \infty $ and $N_{0} = N$ in (\ref{EqPropDescriptionContinuousFunctionals}), one easily shows that the integral in (\ref{eq:SlowGrowingMeasureDualFastDecreasing}) defines a continuous linear functional.   

Let us now prove the converse. Let $T \in C_{FD}'(\mathbb{R}^{d})$. There exists $C >0 $ and $N \in \mathbb{N}$ such that (\ref{EqPropDescriptionContinuousFunctionals}) holds. Let us define the linear functional $(1+|x|^{2})^{-N}T : C_{FD}(\mathbb{R}^{d}) \to \mathbb{C}$ with
\begin{equation}
\label{Definition(1+|x|^{2})^{-N}T}
\langle (1+|x|^{2})^{-N}T , f \rangle := \langle T , (1+|x|^{2})^{-N}f \rangle.
\end{equation}
Since for all $ f \in C_{FD}(\mathbb{R}^{d})$, $ (1+|x|^{2})^{-N}f$ is also in $C_{FD}(\mathbb{R}^{d})$ this functional is well-defined. Considering that $(1+|x|^{2})^{-N}|f| \leq |f|$, it is easy to see that it is continuous. Using (\ref{EqPropDescriptionContinuousFunctionals}) we get
\begin{equation}
\label{Eq(1+|x|^{2})^{-N}TfBounded}
|\langle (1+|x|^{2})^{-N}T , f \rangle | \leq C \| (1+|x|^{2})^{N}(1+|x|^{2})^{-N} f \|_{\infty} = C \| f \|_{\infty},
\end{equation}
for all $f \in C_{FD}(\mathbb{R}^{d})$. In particular, (\ref{Eq(1+|x|^{2})^{-N}TfBounded}) holds for all $f \in C_{c}(\mathbb{R}^{d})$, the space of compactly supported continuous complex functions over $\mathbb{R}^{d}$. Consider $C_{0}(\mathbb{R}^{d})$, the space of continuous complex functions defined over $\mathbb{R}^{d}$ vanishing at infinity, which is a Banach space with the supreme norm. As $C_{c}(\mathbb{R}^{d})$ is a dense subspace of $C_{0}(\mathbb{R}^{d})$, then by extension of bounded linear functionals, we obtain that $(1+|x|^{2})^{-N}T$  is a bounded linear functional over $C_{0}(\mathbb{R}^{d})$, for which (\ref{Eq(1+|x|^{2})^{-N}TfBounded}) holds for every $f \in C_{0}(\mathbb{R}^{d})$. By Riesz's Representation Theorem (see \cite{Rudin1987}, chapter 6) we conclude that $(1+|x|^{2})^{-N}T$ is identified with a unique complex finite measure $\nu$ over $\mathbb{R}^{d}$. By defining $\mu = (1+|x|^{2})^{N}\nu $, we obtain that $\mu$ is a well-defined Radon measure which is in $\mathcal{M}_{SG}(\mathbb{R}^{d})$. Considering that $(1+|x|^{2})^{N}(1+|x|^{2})^{-N}T = T$, it is straightforward that $\langle T , f \rangle = \langle \mu , f \rangle $ for every $f \in C_{FD}(\mathbb{R}^{d})$, and thus $T$ is identified with a unique complex measure in $\mathcal{M}_{SG}(\mathbb{R}^{d}) $. 
$\blacksquare$

\section{Proof of Proposition \ref{Prop:SCEUGenOEM}}
\label{Section:ProofPropositionSCEUGenOEM}

For an arbitrary odd and polynomially bounded continuous function $g_{I}$, let us denote $f_{g_{I}}$ the  associated function defined in Eq. (\ref{eq:SymbolFunctionGenOEM}). Let us first prove the sufficiency of our claim. Let $g_{R}$ satisfy the specified conditions. Let $p : \mathbb{R}^{d} \to \mathbb{R}^{+}_{*}$ be a strictly positive polynomial such that $|g_{R}| \geq 1/p$. When $\beta$ is an odd integer, we get $\cos(\beta\pi/2) = 0$, and it is thus straightforward that $|f_{g_{I}}|^{2} \geq g_{R}^{2} \geq 1/p^{2} $, from which we obtain that $f_{g_{I}}$ satisfies the SCEU independently on $g_{I}$. When $\beta$ is not and odd integer, the choice of the sign of $g_{R}$ is made in order to make that both $\cos(\beta\pi/2) $ and $g_{R}$ have the same sign. Thus, for all $(\xi,\omega) \in \mathbb{R}^{d}\times\mathbb{R} $ we have
\begin{equation}
|f_{g_{I}}(\xi,\omega)|^{2} \geq \left( |\omega|^{\beta}\cos(\frac{\beta\pi}{2}) + g_{R}(\xi)  \right)^{2} = \left( |\omega|^{\beta}|\cos(\frac{\beta\pi}{2})| + |g_{R}|(\xi)  \right)^{2} \geq |g_{R}(\xi)|^{2} \geq \frac{1}{p^{2}(\xi)}.
\end{equation}
Hence, $f_{g_{I}}$ also satisfies the SCEU. Let us now prove the necessity. Suppose that for every $g_{I}$ there exists a strictly positive polynomial $q_{g_{I}} : \mathbb{R}^{d}\times\mathbb{R} \to \mathbb{R}^{+}_{*}$ such that $| f_{g_{I}}| \geq \frac{1}{q_{g_{I}}}$. Then, in particular for $g_{I} = 0$ and evaluating at $\omega = 0$, we get  $ |f_{0}(\xi,0)|^{2} = g_{R}^{2}(\xi) \geq  q_{0}(\xi,0)^{-2} $ from which we obtain that $g_{R}$ satisfies the SCEU. Let $\beta$ be such that $\cos\left(\beta\pi/2\right) < 0$. Suppose there exists $\xi_{1} \in \mathbb{R}^{d}$ such that $g_{R}(\xi_{1}) \geq 0 $. Since $g_{R}$ is continuous and satisfies the SCEU, it follows that $g_{R}>0 $. For every $\xi \in \mathbb{R}^{d}$, if we consider $\omega_{\xi} = (-g_{R}(\xi)/\cos(\beta\pi/2) )^{1/\beta} $ we obtain that $|f_{g_{I}}(\xi,\omega_{\xi})|^{2} = 0 + (-g_{R}(\xi)\tan(\beta\pi/2) + g_{I}(\xi))^{2}$. It is then sufficient to consider a convenient function $g_{I}$ and a point $\xi_{0} \in \mathbb{R}^{d}$ such that $g_{I}(\xi_{0}) = g_{R}(\xi_{0})\tan(\beta\pi/2) $ to obtain $ |f_{g_{I}}(\xi_{0} , \omega_{\xi_{0}} ) |= 0 $, and thus that $f_{g_{I}}$ does not satisfy the SCEU. The contradiction proves that $g_{R}$ must be a negative function. An analogue argument is used to prove that $g_{R}$ must be a positive function when $\beta$ is such that $\cos\left(\beta \pi/2\right) > 0$.$\blacksquare $

\section{Proofs regarding the stochastic heat equation (Ex. \ref{Ex:HeatEquation})}
\label{Section:ProofsHeatEquation}

\subsection{Existence of stationary solutions}
\label{Section:ExistenceSolutionsHeatEquationWhiteNoise}

According to Theorem \ref{TheoremExistAndUniqueStationary}, there exists a stationary solution to the Stochastic Heat Equation with White Noise source term (\ref{eq:SPDEHeatEquationWhiteNoise}) if and only if the spatio-temporal measure $ (\omega^{2}+ a^{2}|\xi|^{4})^{-1}d\xi d\omega $ is in ${\cal M}_{SG}(\mathbb{R}^d \times \mathbb{R})$, i.e. is it is a slow-growing measure. This is the case if the  function $ (\xi,\omega) \mapsto (\omega^{2} + a^{2}|\xi|^{4})^{-1} $ is locally integrable, in which case the slow-growing behavior is provided by the fact that this function is bounded outside a neighborhood around the origin.  It suffices thus to study the integrability over subsets of $\mathbb{R}^{d}\times\mathbb{R}$ of the form $B_{R}^{(d)}(0)\times \left[ - M , M \right]  $ for $R , M > 0$, where $B_{R}^{(d)}(0) \subset \mathbb{R}^{d}$ is the ball of radius $R$ centered in $0$. Using integration with polar coordinates in the spatial domain and the symmetry in the time dimension, we obtain
\begin{equation}
\label{eq:DensityHeatEquationLocalIntegrable}
\int_{B_{R}^{(d)}(0)\times \left[ - M , M \right]} \dfrac{1}{\omega^{2} + a^{2}|\xi|^{4}}  d(\xi,\omega)  = C\int_{0}^{R} \arctan\left(\frac{M}{ar^{2}}\right)r^{d-3}dr
\end{equation}
for some positive constant $C$. Since we have that $\arctan(\frac{M}{aR^{2}}) \leq \arctan(\frac{M}{ar^{2}}) \leq \frac{\pi}{2}$  for all $ r \in \left[ 0 , R \right] $, we conclude that the integral (\ref{eq:DensityHeatEquationLocalIntegrable}) is finite only for $ d > 2 $, from which we get that there exists stationary solutions to the equation (\ref{eq:SPDEHeatEquationWhiteNoise}) only for spatial dimensions $d \geq 3 $. In these cases, solutions would have a functional meaning if the measure $ (\omega^{2}+ a^{2}|\xi|^{4})^{-1}d\xi d\omega $ was finite, which would hold if the limit when $M$ and $R$ go to $\infty$ would exist and was finite. However, by seeing that $\int_{0}^{R}\arctan(\frac{M}{ar^{2}})r^{d-3}dr \geq \arctan(\frac{M}{aR^{2}})\frac{R^{d-2}}{d-2}$, and by letting $M \to \infty$ first and $R \to \infty$ second, one gets that the limit is not finite. Hence, the stationary solutions to (\ref{eq:SPDEHeatEquationWhiteNoise}) only have a meaning as GeRFs and not as Random Functions.

\subsection{Covariance structure}
\label{Section:ObtainingCovarianceHeatEquationWhiteNoised3}

The covariance structure (\ref{eq:CovarianceSPDEHeatEquationWhiteNoise}) is the Fourier transform of the spatio-temporal spectral measure $d\mu_{U}(\xi,\omega) = (2\pi)^{-\frac{d+1}{2}} (\omega^{2}+a^{2}|\xi|^{4})^{-1}d\xi d\omega $ for $d =3$. This measure is not finite. The Fourier transform  $\rho_{U} = \mathscr{F}(\mu_{U})$ is obtained as the limit, in a distributional sense, of continuous functions. Let $R > 0 $ and let denote $\mu_{U}^{R}$  the restriction of the measure $\mu_{U}$ to $B_{R}(0)\times\mathbb{R} \subset \mathbb{R}^{3}\times\mathbb{R}$, i.e. $d\mu_{U}^{R}(\xi , \omega) = (2\pi)^{-\frac{d+1}{2}}(\omega^{2}+a^{2}|\xi|^{4})^{-1}\mathbf{1}_{B_{R}(0)}(\xi)d\xi d\omega  $. This measure is finite, so $\rho_{U}^{R} = \mathscr{F}(\mu_{U}^{R}) $ is a continuous positive-definite function. By the dominated convergence theorem, one gets that for every $\varphi \in \mathcal{S}(\mathbb{R}^{3}\times\mathbb{R})$, $\langle \mu_{U}^{R} , \varphi \rangle \to \langle \mu_{U} , \varphi \rangle $ as $R \to \infty $. Thus, by continuity of the Fourier transform, we have $\rho_{R}^{U} \stackrel{S^{'}}{\to} \rho_{U}$. Let us calculate $\rho_{U}^{R}(h,u) $ for $ (h,u) \in \mathbb{R}^{3}\times\mathbb{R}$.
\begin{align}
\rho_{U}^{R}(h,u) &  = \frac{1}{(2\pi)^{4}}\int_{\mathbb{B}_{R}(0)}\int_{\mathbb{R}}\frac{e^{-iu\omega - i h^{T}\xi}}{\omega^{2} + a^{2}|\xi|^{4}}d\omega d\xi \nonumber  \\
   &  = \frac{1}{(2\pi)^{3}}\frac{1}{2a}\int_{B_{R}(0)}e^{-ih^{T}\xi}\frac{e^{-a|\xi|^{2}|u|}}{a|\xi|^{2}}d\xi \nonumber  \\
   &  = \frac{1}{(2\pi)^{\frac{3}{2}}}\frac{1}{2a}\sqrt{\frac{2}{\pi}}\int_{0}^{R}\frac{J_{\frac{1}{2}}(|h|r)}{\sqrt{r|h|}}e^{-a|u|r^{2}}dr \nonumber \\
   & = \frac{1}{(2\pi)^{2}}\frac{1}{a|h|}\int_{0}^{R}\frac{\sin(|h|r)}{r}e^{-a|u|r^{2}}dr.  \label{eq:RhoURhu} 
\end{align}
Here we have used the expression of the Fourier transform of radial functions. Let us evaluate the limit of $\rho_{R}(h,u)$ when $R \to \infty$ for $|h| \neq 0 \neq |u|$. Consider the function $f_{R} : \mathbb{R}^{+} \to \mathbb{R}$ defined by $f_{R}(\lambda) = \int_{0}^{R}\frac{\sin(\lambda r)}{r} e^{-a|u|r^{2}}dr $ for $\lambda \geq 0 $, with  $f_{R}(0) = 0$. By the dominated convergence theorem, we have $f_{R}'(\lambda) = \int_{0}^{R}\cos(\lambda r) e^{-a|u|r^{2}}dr $. Using the expressions of the Fourier transform of a Gaussian function, one proves that
\begin{equation}
\displaystyle\lim_{R \to \infty} f_{R}'(\lambda) = \sqrt{\frac{\pi}{4a|u|}}e^{-\frac{\lambda^{2}}{4a|u|}}. 
\end{equation}
Using  $f_{R}(\lambda) = \int_{0}^{\lambda}f_{R}'(s)ds $ and the dominated convergence theorem, we get
\begin{equation}
\label{eq:LimfR}
\lim_{R \to \infty} f_{R}(\lambda) = \int_{0}^{\lambda}\sqrt{\frac{\pi}{4a|u|}}e^{-\frac{s^{2}}{4a|u|}}ds =  \frac{\pi}{2}\erf\left(\frac{\lambda}{2\sqrt{a|u|}}\right).
\end{equation}
Using this result in (\ref{eq:RhoURhu}) with $\lambda = |h|$ and $R \to \infty$, we finally obtain the distribution \textit{associated to the function}
\begin{equation}
\label{eq:RhoCovarianceHeatEquationWhiteNoise}
\rho_{U}(h,u)= \dfrac{1}{(2\pi)^{2}} \dfrac{\pi}{2a|h|}\erf\left( \frac{|h|}{2\sqrt{a|u|}} \right). 
\end{equation}
which is the expression in (\ref{eq:CovarianceSPDEHeatEquationWhiteNoise}). 

It is worth emphasizing that this expression is only valid in a distributional sense. The distribution $\rho_{U}$ is only meaningful when applied to test functions, satisfying 
\begin{equation}
\langle\rho_{U} , \psi \rangle = \displaystyle\lim_{R \to \infty} \langle \rho_{U}^{R} , \psi \rangle, \quad  \forall \psi \in \mathcal{S}(\mathbb{R}^{3}\times\mathbb{R}). 
\end{equation}
The expression \textit{associated to the function (\ref{eq:RhoCovarianceHeatEquationWhiteNoise})} refers to the fact that for every test function $\psi$ such that its support does not contain the origin, we have
\begin{equation}
\langle \rho_{U} , \psi \rangle = \int_{\mathbb{R}^{3}\times\mathbb{R}} \frac{1}{(2\pi)^{2}} \frac{\pi}{2a|h|}\erf\left( \frac{|h|}{2\sqrt{a|u|}} \right) \psi(h,u)dhdu.
\end{equation}

\end{appendices}
\bibliographystyle{vancouver}

\end{document}